\documentclass[a4paper,10pt]{article}
\usepackage[cp1251]{inputenc}
\usepackage[english]{babel}

\usepackage{amssymb}\usepackage{amsmath}
\usepackage{graphicx}
\def\bl{\rule[-1mm]{2.4mm}{2.4mm}}

\def\be{\begin{equation}}
\def\ee{\end{equation}}

\def\ord{{\rm ord}}
\def\bl{\rule[-1mm]{2.4mm}{2.4mm}}
\def\vert{{\rule{.5mm}{2.mm}}}
\def\hor{{\rule{2.mm}{.5mm}}}

\newtheorem{thrm}{\bf Theorem}
\newtheorem{lmm}{\bf Lemma}
\newtheorem{dfn}{\bf Definition}
\newtheorem{rmk}{\bf Remark}
\newtheorem{xmpl}{\bf Example}
\numberwithin{rmk}{section}
\numberwithin{lmm}{section}
\numberwithin{xmpl}{section}

\begin{document}

\title {Combinatorial analysis of the period mapping:\\ topology of 2D fibers}
\author{\copyright 2017 ~~~~A.B.Bogatyrev
\thanks{Supported by RFBR grants 16-01-00568, and RAS Program
"Modern problems of theoretical mathematics"}}
\date{}
\maketitle
{\bf Abstract:} {\small
We study the period mapping from the moduli space of real hyperelliptic
curves with marked point on an oriented oval to the euclidean space. The mapping
arises in the analysis of Chebyshev construction used in the 
constrained optimization of the uniform norm of polynomials and rational
functions. The decomposition of the moduli space into polyhedra labeled by planar graphs 
allows to investigate the global topology of low
dimensional fibers of the period mapping.}

\section{Introduction}
The period mappings are defined on the moduli spaces of curves with marked points on them and return the values of periods of 
uniquely associated to such curves differentials with singularities in the marked points.
They appear in different settings \cite{AVG, Kr, KrGr, KuiMo},
mostly for the study of the geometry of the moduli spaces of curves. For instance in \cite{Kr}, the authors try to construct 
the subvarieties of the moduli space as the closed leaves of the foliation locally defined by the period mapping.  
The interrelation of the "integral" leaves of sertain foliations with the spectral curves of elliptic Calogero-Moser 
systems is shown in \cite{Kri1}. In \cite{KuiMo} the authors use period mapping to study the asymptotics of orthogonal polynomials.
This particular research was motivated by uniform rational approximation problems.

A universal phenomenon known as "Chebyshev Alternation Principle" in mathematical community or
"Equiripple Property"  in the community of electrical engineers says that solutions to certain 
uniform norm optimization problems for polynomials look like 
 waves of constant amplitude. The most known function with this property is the (co)sine, so there is possibly no wonder that
the degree $n$ solutions $P_n(x)$ may be effectively evaluated by finding parameters in the following Chebyshev representation (Ansatz)
\cite{B0} and \cite{BBook} Chap.2:
$$
P_n(x)=\cos(ni\int_{(e,0)}^{(x,w)}d\eta_M),
$$
where $d\eta_M$ is a suitably normalized third kind abelian differential on a hyperelliptic curve $M$.
For the above formula to give a single valued function on a Riemann sphere, all the periods of the abelian integral 
should be multiples of the periods of cosine, in other words lie in the lattice $2\pi in^{-1}\mathbb{Z}$.
Deformation of the initial optimization problem brings us to the deformation of the curve $M$ and the induced deformation 
of the differential $d\eta_M$ such that all its periods are conserved since one cannot continuously jump 
from one point of the lattice to the other. The above representation of polynomials in terms of Riemann surfaces also gives the solution to 
Pell-Abel equation which in turn is related to Poncelet problem, elliptic billiards \cite{VIDr} and boundary value problem for string equation \cite{BZ}.  

The brief content of this paper is as follows.
In Sect. 2 we introduce the moduli spaces of real hyperelliptic curves, a model of their universal covering -- the labyrinth space --
and define the period mapping from the labyrinth to the euclidean space. The main object of our research is the  
global topology of the fibers of the latter map. In Sect. 3 we manufacture the main tool of our investigation, the description of 
curves in terms of weighted flat graphs composed of the segments of two orthogonal foliations intrinsically associated to the curve. In Sect. 4 we study the decomposition of the moduli spaces into polyhedra spanned by the weights of the admissible graphs.
Sect. 5 presents the period mapping as linear one  in terms of local coordinates of each polyhedron.
In Sect. 6 we construct a PL model for fibers of the periods map attaching one to the other polyhedra being the sections of the polyhedra 
from Sect. 4 by linear spaces. In particular, we conclude that the topology of 2D fibers is trivial.   

{\bf Acknowledgements}. The paper was started when the author hold a visiting position at the University Bordeaux-I (LABRI) which provided excellent conditions for work and interdisciplinary communication. My special thanks are to Sasha Zvonkine who instigated my interest in this topic after it nearly extinguished 
and who gently pushed me toward writing this manuscript. 

\section{Setting the problem}
\subsection{Moduli spaces}
Let $\cal H$ be the moduli space of smooth real hyperelliptic curves $M$ with one marked point $"\infty"$ (in what follows will be used without quotes) on an oriented real oval.
We require this marked point being not fixed by the hyperelliptic involution $J$ acting on each curve. Fixing the
genus $g=0,1,\dots$ of a curve and the number $k=0,1,\dots,g+1$ of its real ovals, we arrive to the decomposition of the moduli space into 
connected components with given values of the two topological invariants:
$$
{\cal H}=\sqcup_{g,k}{\cal H}_g^k.
$$

Any element of $\cal H$ admits an affine model:  
\be
\label{M}
M=M({\sf E}):= \{(x,w)\in\mathbb{C}^2:\quad w^2=\prod_{s=1}^{2g+2}(x-e_s)\}, 
\qquad e_s\in {\sf E},
\ee
with branching set $\sf E=\bar{E}$ symmetric with respect to the real axis. 
The hyperelliptic and anticonformal involutions of the curve we define as
$J(x,w):=(x,-w)$ and $\bar{J}(x,w):=(\bar{x},\bar{w})$ respectively.
The marked point on the real oval is the point corresponding to $(x,w)=(+\infty,+\infty)$ in the natural two-point compactification of
the curve (\ref{M}). 

An element in the component ${\cal H}_g^k$ of the total moduli space is represented by a set $\sf E^+$ of $g-k+1$ distinct points in the open upper half plane $\mathbb{H}$ and $2k$ distinct points on the real axis $\mathbb{R}$ modulo translations and dilations of the set. In this paper we only consider the case $k>0$, and the set 
$\sf E^+$ may be normalized so that its two extreme real points are $\pm1$. In this notation the branching set $\sf E$ of the curve (\ref{M})
is the union of the set $\sf E^+$ and its complex conjugate. 

\begin{lmm}\cite{B2},\cite{BBook} {\rm ~Chap. 3.}~~~
The space ${\cal H}_g^k$ is a smooth real manifold of dimension $2g$ and its fundamental group is the group  $Br_{g-k+1}$ of braids on $g-k+1$ strands. 
\end{lmm}

\subsection{Local period mapping}
On each curve $M$ from the moduli space $\cal H$ there is a unique third kind abelian differential $d\eta_M$
with just two simple poles: the marked point $\infty$ and its involution $J\infty$, with residues $-1$ and $+1$ respectively and  purely imaginary periods \cite{BBook}, \S 2.1.3. For the algebraic model (\ref{M}) of the curve $M$ the differential takes the form:
\be
\label{deta}
d\eta_M=(x^g+\dots)w^{-1}dx,
\ee
with dots standing for a polynomial of degree not greater than $g-1$. One can check that normalization conditions of the differential imply that the latter is \emph{real}, that is $\bar{J}d\eta_M=\overline{d\eta_M}$. In other words, the polynomial in (\ref{deta})
has real coefficients. An important consequence of this fact is this: half of periods of this differential vanish.

Indeed, the anticonformal involution $\bar{J}$ naturally acts on the space of cycles $H_1(M,\mathbb{R})$ and splits it into a sum of two eigenspaces 
$H_1(M,\mathbb{R})^\pm$ corresponding to eigenvalues $\pm1$ of the operator. Cycles invariant under the reflection 
$\bar{J}$ we call \emph{even} and those changing the sign we call \emph{odd}. The periods of any real differential $d\eta$ over even/odd cycles are respectively real/ purely imaginary:
$$
\int_C d\eta=\pm\int_{\bar{J}C} d\eta=
\pm\int_C \bar{J}d\eta=\pm\int_C \overline{d\eta}=
\pm\overline{\int_C d\eta}.
$$
In particular, we see that the integrals of the distinguished differential $d\eta_M$
over even cycles vanish as they should be real and imaginary simultaneously.

One can locally define the period mapping in the neighbourhood of a distinguished point $M_0$ of the moduli space 
as 
\be\label{PM}
\langle
\Pi(M)|C
\rangle
:=i\int_C d\eta_M,  
\qquad C\in H_1^-(M_0,\mathbb{Z}),
\ee
here the integration of the differential living on a curve $M$ over the cycles on the distinguished curve $M_0$ is possible due to 
Gauss-Manin connection in the homological bundle over the moduli space \cite{Vas}, \cite{BBook} Chap.5. Globally this mapping is not well defined due to nontrivial holonomy of the connection: braids entangle the cycles on the curve. We study this effect after a little while.

\subsection{Universal covering of the moduli space}
The usual way to correctly define a multivalued mapping is to lift it to the universal covering of the source space. 
The space $\tilde{\cal H}_g^k$ -- the universal covering of ${\cal H}_g^k$ -- has several models described in \cite{B1}, \cite{BBook}, Chap.3.
For instance, it may be represented by Teichm\"uller space of a $(g-k+1)$ -times punctured disc with marked points on the boundary. In particular,  the universal covering is homeomorphic to the euclidean space $\mathbb{R}^{2g}$. Here we use another model of this space, namely the space of labyrinths.

\begin{figure}
\begin{picture}(110,40)
\put(0,5){
\begin{picture}(50,40)
\put(0,0){\line(1,0){45}}
\multiput(5,0)(5,0){4}{\circle*{1}}
\thicklines

\put(40,10){\circle*{1}}
\put(35,15){\circle*{1}}
\put(10,10){\circle*{1}}
\put(21,13){$\mathbb{H}$}
\end{picture}}

\put(50,5){
\begin{picture}(60,40)
\put(0,0){\line(1,0){55}}
\multiput(5,0)(5,0){4}{\circle*{1}}
\put(21,13){$\mathbb{H}$}
\put(6,-4){$\Lambda_0$}
\put(16,-4){$\Lambda_1$}

\thicklines
\multiput(5,0.3)(10,0){2}{\line(1,0){5}}
\put(25,0){\oval(10,10)[tr]}
\put(25,15){\oval(20,20)[l]}
\put(25,15){\oval(40,20)[tr]}
\put(40,15){\oval(10,10)[br]}
\put(40,10){\circle*{1}}
\put(29,-4){$\Lambda_2$}

\put(35,0){\oval(10,10)[tr]}
\put(35,10){\oval(10,10)[l]}
\put(35,15){\circle*{1}}
\put(39,-4){$\Lambda_3$}

\put(30,10){\oval(40,40)[t]}
\put(50,0){\line(0,1){10}}
\put(10,10){\circle*{1}}
\put(49,-4){$\Lambda_4$}
\end{picture}}
\end{picture}
\includegraphics[scale=.45, trim = 0cm -1cm 0cm 0cm, clip]{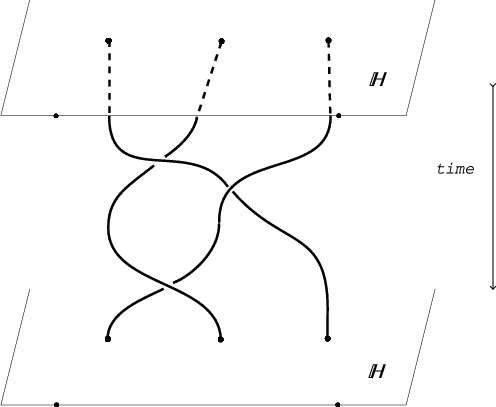}
\caption{\small  A point  ${\sf E}^+\in{\cal H}^2_4$ (left) is lifted to the universal cover 
by choosing the labyrinth that escorts it (middle). A loop in the moduli space ${\cal H}^1_3$ is a braid on three strings (right)}
\label{Labrnth}
\end{figure}
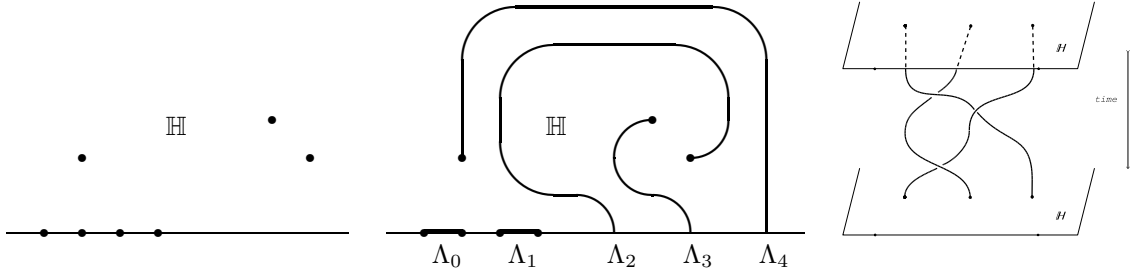

A labyrinth $\Lambda=(\Lambda_0,\dots,\Lambda_g)$ attached to a point $\sf E^+$ of the moduli space (see Fig. \ref{Labrnth}) is a set of $g+1$ disjoint  arcs  in the closed upper half plane $\hat{\mathbb{H}}:=\mathbb{H}\cup\mathbb{R}\cup\infty$. First 
$k$ of them are real segments $\Lambda_0,\dots,\Lambda_{k-1}$ pairwise connecting real points of $\sf E^+$. 
The remaining $g-k+1$ arcs $\Lambda_k$, $\Lambda_{k+1}$, $\dots\Lambda_g$ have more freedom: they start at the complex points of the branching set $\sf E^+$ and meet 
the real axis to the right of the largest real point of the set.
All the arcs of a labyrinth are ordered by their intersection with the real axis.
Two labyrinths $\Lambda$  obtained by the isotopy of the upper half plane fixed at the points of $\sf E^+$ are declared to be equal. 

\begin{dfn}
The space of labyrinths ${\cal L}_g^k$ is the set of the equivalence classes of labyrinths $\Lambda$  attached to the points $\sf E^+$ running through the moduli space component ${\cal H}_g^k$. The total labyrinth space is the disjoint union of the connected components:
$$
{\cal L}:=\mathop{\sqcup}\limits_{g,k}{\cal L}_g^k,
$$
a useful notion if we do not want to specify the values of the topological invariants $g,k$ of a real curve.
\end{dfn}

Speaking informally, the universal covering of a space is the history of motion from some marked point told for each point of this space. The variable part $\Lambda_k,\dots,\Lambda_g$ of a labyrinth has the meaning of the traces of the motion of branch points in the upper half plane.

The braid group $Br_{g-k+1}$ realized as the mapping class group \cite{Bir} of the punctured half plane $\mathbb{H}\setminus{\sf E}$ naturally acts on the space of labyrinths ${\cal L}_g^k$  and corresponds to cover transformations. Projection from the universal cover to the (component of the)  moduli space consists in wiping out the labyrinth leaving the branch points only.

\subsection{Distinguished basis in odd homologies}
Any labyrinth $\Lambda$ gives us a distinguished basis in the lattice of odd cycles on the surface $M$.
Take a closed upper half plane $\hat{\mathbb{H}}$ with removed labyrinth $\Lambda$.
This set is simply connected and avoids the branch points of the two sheeted covering $M\to M/J=\mathbb{C}P^1$.
Therefore, the set may be uniquely lifted to the surface by the requirement that infinity of the cut plane is mapped to the marked point 
$"\infty"$ of the surface. The union of $\hat{\mathbb{H}}\setminus\Lambda$ embedded in the surface $M$ and its 
reflection $\bar{J}(\hat{\mathbb{H}}\setminus\Lambda)$ is a subsurface of $M$ with $g+1$ boundary components which we endow 
with the standard boundary orientation\footnote{Opposite to the notations of \cite{BBook}, we've also changed the sign in the 
definition (\ref{PM}) of the periods map.}. Let us call those cycles $C_0$, $C_1$, $\dots,C_g$, enumeration is inherited from that of the 
labyrinth components.

\begin{lmm}\cite{B0, B2},\cite{BBook}, {\rm Chap.2~~~} 
\label{OddBasis}
1)The associated cycles  $C_0$, $C_1$, $\dots,C_g$ of a labyrinth  make up a basis in the lattice of 
odd cycles on the surface $M({\sf E})$ punctured at two points $\infty$ and $J\infty$. \\
2) Same cycles make up horizontal (with respect to Gauss-Manin connection) sections of homology bundle
over universal cover of the moduli space.
\end{lmm}
 
\begin{figure}

\begin{picture}(170,40)
\put(10,5){
\begin{picture}(60,40)
\put(0,0){\line(1,0){55}}
\multiput(15,35)(25,0){2}{\circle*{1}}
\thicklines
\multiput(15,0)(25,0){2}{\line(0,1){35}}
\thinlines
\multiput(15,0)(25,0){2}{\oval(6,76)[t]}
\multiput(18,20)(25,0){2}{\vector(0,-1){5}}
\multiput(12,20)(25,0){2}{\vector(0,1){5}}
\put(14,-7){$C_s$} 
\put(39,-7){$C_{s+1}$} 
\put(26,5){$\mathbb{H}$} 
\put(63,-1){$\mathbb{R}$} 
\end{picture}}

\put(80,5){
\begin{picture}(60,40)
\put(0,0){\line(1,0){75}}
\multiput(15,35)(25,0){2}{\circle*{1}}
\thicklines
\put(27.5,10){\oval(25,20)[tl]}
\put(27.5,35){\oval(25,30)[rb]}
\put(15,0){\line(0,1){10}}

\put(35,35){\oval(40,20)[t]}
\put(55,0){\line(0,1){35}}

\thinlines
\put(58,0){\line(0,1){38}}
\qbezier(58,38)(58,48)(48,48)

\put(58,18){\vector(0,-1){5}}

\put(49,48){\line(-1,0){27}}
\qbezier(22,48)(12,48)(12,38)
\put(12,39){\line(0,-1){4}}
\qbezier(12,35)(12,32)(15,32)
\qbezier(15,32)(18,32)(18,35)
\put(18,35){\line(0,1){2}}
\qbezier(18,37)(18,42)(23,42)
\put(23,42){\line(1,0){24}}
\qbezier(47,42)(52,42)(52,37)
\put(52,37){\line(0,-1){37}}
\put(52,18){\vector(0,1){5}}

\put(12,0){\vector(0,1){12}}
\qbezier(12,12)(12,23)(23,23)

\put(23,23){\line(1,0){9}}
\qbezier(32,23)(37,23)(37,28)

\put(37,28){\line(0,1){7}}
\qbezier(37,35)(37,38)(40,38)
\qbezier(40,38)(43,38)(43,35)

\put(43,35){\vector(0,-1){10}}
\qbezier(43,28)(43,17)(32,17)

\put(32,17){\line(-1,0){9}}
\qbezier(23,17)(18,17)(18,12)
\put(18,12){\line(0,-1){12}}

\put(35,5){$\mathbb{H}$} 
\put(14,-7){$C_s$} 
\put(61,-7){$C_{s+1}$} 

\end{picture}}
\end{picture}

\caption{\small  An elementary braid $\beta_{s-k+1}$, $s=k,g-1$, twists the labyrinth and changes the distinguished basis of odd cycles}
\label{Burau}
\end{figure}
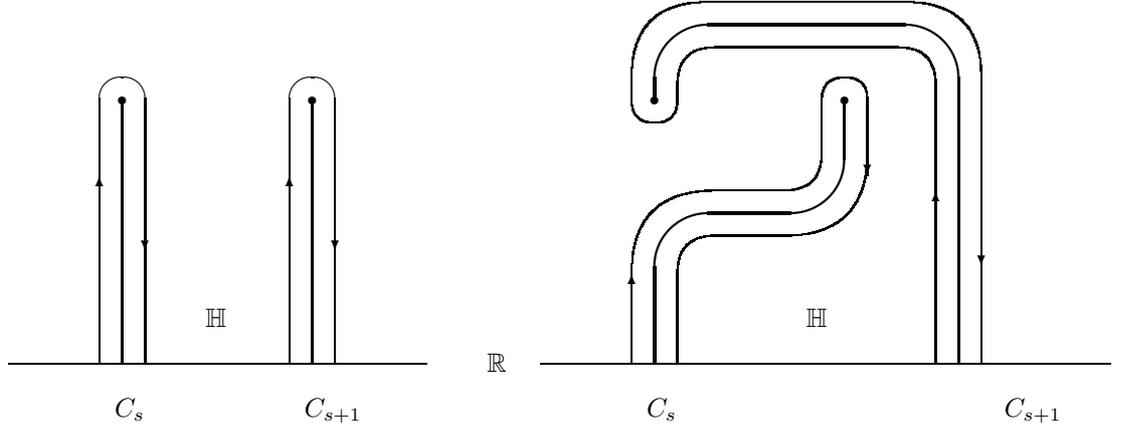

The distinguished basis of homologies changes under the group of cover transformations $Br_{g-k+1}$ of
the universal cover $\tilde{\cal H}^k_g$ in a predictable way (holonomy of GM connection).
An elementary braid $\beta_{s-k+1}\in Br_{g-k+1}$, $s=k,\dots,g-1$, is represented by a counterclockwise half-twist along a contour only intersecting the neighboring 
arcs $\Lambda_{s}$ and $\Lambda_{s+1}$ of the labyrinth. One can check \cite{BBook,B2} that all the cycles of the distinguished basis of homologies remain intact, but two of them (see Fig. \ref{Burau} ):
\be
\beta_{s-k+1}:\qquad (C_{s}, C_{s+1})^t\to (- C_{s+1}, 2C_{s+1}+ C_{s})^t.
\label{BuRep}
\ee
This matrix representation of the braid group is known as (a particular case of) Burau representation \cite{Bir}.

\subsection{Global period mapping}
\begin{dfn}
The periods map (\ref{PM}) evaluated at the distinguished basis of odd cycles gives a well-defined 
global mapping 
$$
\Pi(\Lambda): \quad {\cal L}_g^k\to\mathbb{R}^{g+1}.
$$
 
\end{dfn}

The image of this map lies in a codimension one euclidean subspace: the sum of the distinguished cycles
is homologous to the circle encompassing the pole of the differential $d\eta_M$, therefore
the value of the sum of the periods map components is always $2\pi$.

\begin{thrm} \cite{B1},\cite{BBook}, {\rm Chap. 5~~~}\\
1) The period mapping $\Pi(\Lambda): \quad {\cal L}_g^k\to\mathbb{R}^g$ has full rank.\\
2) period mapping is equivariant with respect to the action of braids:
$$
\beta\Pi=\Pi\beta, 
\qquad \beta\in Br_{g-k+1}.
$$
Here braids act on the universal covering of the moduli space as cover transformations and on the euclidean space 
by Burau representation (\ref{BuRep}).
\label{PMfacts}
\end{thrm}
First statement of the theorem was proved in a more general setting in \cite{Kr,KrGr}.
The range of this periods map was explicitly calculated in \cite{B2},\cite{BBook}, \S 5.3.

\subsection{Statement of the main result}
In this paper we study the topology of fibers of the periods map $\Pi$, i.e. of the inverse images of points from the range of  the mapping. It follows from the above Theorem \ref{PMfacts} that the fibers are smooth $g$-dimensional submanifolds of ${\cal L}_g^k\cong\mathbb{R}^{2g}$. In the next section we introduce a machinery that reduces the problem under investigation to certain combinatorial calculations. In the concluding sections we perform the calculations for the two-dimensional fibers and arrive at the main

\begin{thrm}
\label{MainTh}
Any fiber of the periods map defined in the space ${\cal L}_2^k$, $k=1,2,3$, is a cell.
\end{thrm}

\begin{rmk}
Similar calculations were performed by the author for the 3D fibers too and a new effect was discovered. Fibers may be disconnected, 
however each component of a fiber remains a cell (= topological space homeomorphic to a euclidean one). 
This observation leads us to a \\
{\bf Conjecture:} \cite{B3} Components of the fibers of the period mapping are always cells.
\end{rmk}

\section{Pictorial representation of curves}
Here we work out the main tool of our investigation. We need a convenient description of curves which allows us to effectively reconstruct the period mapping. The construction is briefly described  below, more details may be found in \cite{B2, BBook}.

The idea to represent algebraic curves by (weighted) graphs is not new.
Possibly, the seminal input is due to A.Grothendieck with his Dessins d'Enfants.
Once appeared in math physics, ribbon graphs today make up a flourishing industry established in 
the works of M.Kontsevich, R.Penner, L.Chekhov, V.Fock \cite{Kon,Pen,Fock} to name a few. 
Similar ideas were used by M.Bertola in his work on Boutroux curves \cite{Bert}, see also 
\cite{FV, Sol} . We need a graphical description of curves which allows us to effectively reconstruct the period mapping.
The construction is briefly described  below, more details may be found in \cite{B2}, \cite{BBook}, \S 4.1.

\subsection{Global width function}
Suppose $M({\sf E})\in{\cal H}$ and $d\eta_M$ is the 3rd kind differential  associated with the curve $M$ as above.
One immediately checks that the normalization conditions of $d\eta_M$ force the width function 
\be
\label{W}
W(x):=|Re\int_{(e,0)}^{(x,w)}d\eta_M|,
\qquad x\in\mathbb{C},
\qquad e\in {\sf E},
\ee
to obey the following properties:
\begin{itemize}
\item $W$ is single valued on the plane, 
\item $W$ is harmonic outside its zero set $\Gamma_\vert:=\{x\in\mathbb{C}: W(x)=0\}$, 
\item $W$ has a logarithmic pole at infinity,
\item $W$ vanishes at each branchpoint $e_s\in\sf E$\\
(and hence the definition (\ref{W}) is independent  of the choice of the branchpoint $e$).
\end{itemize}

We only comment on the last property. Since $d\eta_M$ is odd with respect to the hypereliptic involution of $M$,
$W(e_s)$ is equal to one half of the absolute value of the real part of some period of the differential.
Normalization implies that all its periods are purely imaginary.

\subsection{Construction of the graph $\Gamma(M)$.}

To any curve $M\in{\cal H}$ we associate a weighted planar graph $\Gamma=\Gamma(M)$ composed of the finite number of segments of vertical and horizontal foliations \cite{Str} of the quadratic differential $(d\eta_M)^2$ descended to the Riemann sphere. The graph $\Gamma(M)$ is a union of the 'vertical' subgraph $\Gamma_\vert$ and the 'horizontal' subgraph $\Gamma_\hor$ -- see left panel at Fig.\ref{ExtGraph} for an admissible graph. Sets of edges and vertexes of a graph $\Gamma$ are denoted as $Sk_1(\Gamma)$ and $Sk_0(\Gamma)$ respectively.

\begin{dfn}
\begin{itemize} 

\item  VERTICAL EDGES  $:=Sk_1(\Gamma_\vert)$  are arcs of the  zero set of $W(x)$; they are segments of the vertical foliation $d\eta_M^2<0$ and are not oriented.  

\item HORIZONTAL EDGES $:=Sk_1(\Gamma_\hor)$  are all segments of the horizontal foliation $(d\eta_M)^2>0$ (or steepest descent lines for $W(x)$) connecting the finite critical points of the foliation either to other such points or -- as a rule -- to the zero set of $W$. Horizontal edges are oriented with respect to the growth of $W(x)$. 

\item WEIGHTS: Each edge, no matter what type is it, is equipped with its length in the metric $ds=|d\eta_M|$ of quadratic differential.

\item  VERTEXES $:=Sk_0(\Gamma)$ of the graph $\Gamma$ comprise all finite points of the divisor of the quadratic differential $({d\eta_M})^2$
considered on the plane as well as points in $\Gamma_\vert\cap\Gamma_\hor$ -- projections of the saddle points of $W$ to its zero set 
along the horizontal leaves. 
\end{itemize}
\end{dfn}

Instead of assigning lengths to the horizontal edges, it is more convenient to keep the values of the width function
$W(x)$ at all vertexes of the graph: the length of the oriented edge thus is the increment of the width function along it.

\begin{rmk}
The most important property of the description of curves in terms of their graphs $\Gamma$ is the following.
The periods of the differential $d\eta_M$ are integer linear combinations of the lengths of the vertical edges.
Indeed, the weight of any edge is the absolute value of the abelian integral taken along this edge.
Given a cycle on $M$, one just has to properly collapse it to the graph $\Gamma$ to get the period of $d\eta_M$
along this cycle. We consider this in greater detail in Sect. \ref{Sec5}.
\end{rmk}

\begin{rmk}
The topological invariants $g,k$ of the curve $M$ may be reconstructed from the combinatorics of the graph $\Gamma(M)$.
Indeed, one immediately checks that the multiplicity of a vertex $V$ of the graph in the divisor of the quadratic differential 
$(d\eta_M)^2$ equals to 
$$
\ord (V):= d_\vert(V)+2d_{in}(V)-2,
$$
where $d_\vert$ is the degree of the vertex with respect to the vertical edges and $d_{in}$ is the number of incoming horizontal edges. Branchpoints of the curve $M$ correspond to the vertexes with odd value of $\ord(V)$ (or odd value of $d_\vert(V)$ which is the same).  The number $2k$ of branchpoints on the real axis determines the number of real ovals of the curve; the total number $2g+2$ of branchpoints is related to the genus.
\end{rmk}
\begin{dfn}
We call a vertex $V\in Sk_0(\Gamma)$ a branchpoint, if the value  $\ord(V)$ is odd.  
\end{dfn}

\subsection{Axiomatic description of graphs}
We distinguish between the geometric graph $\Gamma(M)$ drawn on the plane,
the class of weighted planar graphs $\{\Gamma\}$ modulo isotopies of the plane respecting the complex conjugation
and the class of the topological planar graphs $[\Gamma]$ without weights.
We are going to describe all admissible types of graphs $\{\Gamma\}$ in axiomatic way
so that they could be listed by an automaton.
It turns out that there are just three topological restrictions on the graph $[\Gamma]$ 
and two more normalization conditions on its weights. All of them are listed in the following 

\begin{lmm}\cite{B2},\cite{BBook}, \S 4.1.4\\
(T1) $\Gamma$ is a tree with real symmetry axis.

(T2) Horizontal edges leaving the same vertex are separated by a vertical or an incoming edge, 
in particular there are no hanging horizontal vertexes like~~ 
\begin{picture}(15,5)
\put(1,0){\circle*{1}}
\put(1,0){\vector(1,0){10}}
\end{picture}
 
(T3) If $\ord (V)=0$ then $V\in\Gamma_\hor\cap\Gamma_\vert$. 

(W1) Width function increases along oriented edges and $W(V)=0$ if $V$ lies on the vertical part of the graph. 

(W2) The weights of vertical edges are positive and their total sum is $\pi$. 
\label{five}
\end{lmm}
{\bf Proof sketch}
We say a few words about properties T1, T2, W2. The rest follow from the definition of the graph $\Gamma$.

(T1). Suppose there is a cycle in the graph. Let us calculate the Dirichlet integral for the width function in the domain $\Omega$
bounded by the cycle by means of Green's formula:
$$
\int_\Omega |grad~W(x)|^2 d\Omega = 
\int_{\partial\Omega} W(x)\frac{\partial W}{\partial n} ~ds. 
$$
Function $W$ vanishes on the vertical parts of the boundary while its normal derivative vanishes at the horizontal parts of ${\partial\Omega}$, 
therefore $W$ is constant. Now suppose the graph has several components. Sum up the value 
$\ord(V)$ over all the vertexes. We get the number $2\sharp\{vertical~edges\}+2\sharp\{horizontal~edges\}-2\sharp\{vertices\}=
-2\sharp\{trees~in~the~forest\}$. This is equal to the degree of the divisor of the quadratic differential $d\eta_M^2$
(i.e. $-4$) plus the order of its pole at infinity (i.e. $2$). Hence, the graph $\Gamma$ is a tree.

(T2). Suppose, $W(V)>0$ at the vertex $V$ of the graph.  This is a saddle point of the width function, the meeting point of 
several alternating ``ridges'' and ``valleys''. A horizontal edge comes into $V$ from each valley. The outgoing edge (if any) goes along the ridge, so any two of them are separated. Same is true for $W(V)=0$ with vertical edges coming from each ``valley''.

(W2). The integral of $d\eta_M$ along the boundary of the plane cut along $\Gamma_\vert$ equals $2i$ times the sum of the weights of all vertical edges.   
The integration path may be contracted to the path encompassing the pole at infinity, hence the integral is $2\pi i$. 
~~~\bl

It turns out that there are no any further restrictions neither on topology nor on weights of graphs:

\begin{thrm}\cite{B2},\cite{BBook}, \S 4.1.5
Each weighted graph 
$\{\Gamma\}$ satisfying five properties of Lemma \ref{five}
stems from a unique curve $M=M\{\Gamma\}\in{\cal H}$.
\label{CfromG}
\end{thrm}

{\bf Proof sketch.} The Riemann surface $M$ may be glued from a finite number of stripes
in a way determined by combinatorics and weights of the graph. Below we briefly describe the procedure.

\begin{figure}
\begin{picture}(80,40)
\unitlength=.8mm
\thinlines
\put(10,20.5){\line(1,0){15}}
\put(10,19.5){\line(1,0){15}}
\multiput(40,35.5)(0,-1){2}{\line(1,0){30}}
\multiput(40,5.5)(0,-1){2}{\line(1,0){30}}
\thicklines
\put(25,20){\vector(1,0){15}}
\put(55,20){\vector(-1,0){15}}
\put(55,35){\vector(0,-1){15}}
\put(55,5){\vector(0,1){15}}

\multiput(10,20)(15,0){4}{\circle*{2}}
\multiput(40,35)(15,0){3}{\circle*{2}}
\multiput(40,5)(15,0){3}{\circle*{2}}
\put(33,5){$\small-1$}
\put(72,5){$\small-1$}
\put(54,0){$\small0$}
\put(57,15){$\small2$}
\put(39,15){$\small2$}
\put(6,15){$\small-1$}
\put(21,15){$\small-1$}

\multiput(10,20)(-1.5,0){7}{.}
\multiput(55,20)(1.5,0){7}{.}
\put(13,23){$H_3$}
\put(44,38){$H_2$}
\put(58,38){$H_1$}
\put(36,23){$W_2$}
\put(57,23){$W_1$}
\end{picture}
\hfill
\begin{picture}(80,40)
\unitlength=.8mm
\thinlines
\put(10,20.5){\line(1,0){15}}
\put(10,19.5){\line(1,0){15}}
\multiput(40,35.5)(0,-1){2}{\line(1,0){30}}
\multiput(40,5.5)(0,-1){2}{\line(1,0){30}}
\thicklines
\put(25,20){\vector(1,0){15}}
\put(55,20){\vector(-1,0){15}}
\put(55,35){\vector(0,-1){15}}
\put(55,5){\vector(0,1){15}}

\multiput(10,20)(15,0){4}{\circle*{2}}
\multiput(40,35)(15,0){3}{\circle*{2}}
\multiput(40,5)(15,0){3}{\circle*{2}}

\thinlines
\put(40,40){\oval(10,10)[bl]}
\put(35,40){\vector(0,1){10}}
\put(70,40){\oval(10,10)[br]}
\put(75,40){\vector(0,1){10}}
\put(55,35){\vector(0,1){15}}

\put(40,0){\oval(10,10)[tl]}
\put(35,0){\vector(0,-1){10}}
\put(70,0){\oval(10,10)[tr]}
\put(75,0){\vector(0,-1){10}}
\put(55,5){\vector(0,-1){15}}
\put(55,20){\vector(1,0){30}}
\put(10,20){\vector(-1,0){10}}
\put(25,20){\oval(30,20)[r]}
\put(25,30){\vector(-1,0){25}}
\put(25,10){\vector(-1,0){25}}
\end{picture}
\vspace{1mm}
\caption{\small Graph $\Gamma(M)$ for a curve $M\in{\cal H}_2^1$ (left) and it's extension $Ext~\Gamma$ (right). Double lines/arrows are vertical/horizontal edges; dotted line is the real (mirror symmetry) axis; numbers at the vertexes stand for the order of the quadratic differential $(d\eta^2_M)$ at those points}
\label{ExtGraph}
\end{figure}
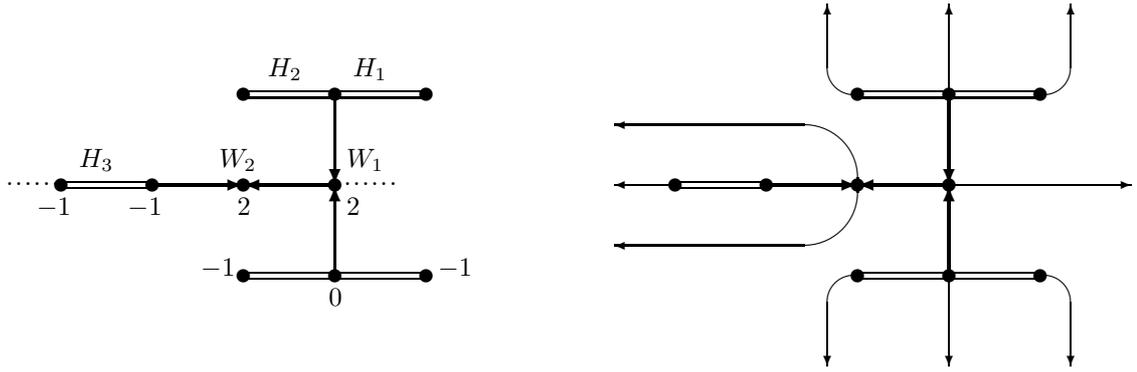

Given a planar graph satisfying all the requirements of Lemma \ref{five}, we extend it by drawing $d_\vert(V)-d_{out}(V)+d_{in}(V)\ge0$ outgoing 
horizontal arcs which connect each vertex $V$ to infinity and are disjoint except possibly at their endpoints. For each vertex of the 
extended graph $Ext~\Gamma$
we require that all the outgoing edges, old and new, alternate with the incident edges of other types: incoming or vertical. 
In particular, the property (T2) of the Lemma \ref{five} is kept for $Ext~\Gamma$. 
Up to isotopy of the plane, the extension of the graph is unique  since the original graph is a tree -- see the right 
picture in Fig. \ref{ExtGraph}. 

From the topological viewpoint all the components of the complement to the extended graph in the plane have the same structure -- 
see right panel at Fig. \ref{ExtGraph}. They are 2-cells bounded by exactly one vertical edge with attached to the endpoints of the 
latter two chains of horizontal edges all pointing in the same direction and eventually meeting at infinity. One easily checks 
that for the graph $\Gamma(M)$ generated by an element of the moduli space, and its extension drawn by the horizontal trajectories, 
the (suitably chosen branch of) abelian integral $\eta(x)=\int^x d\eta_M$ maps each of the above cells to a horizontal half strip of the height equal to the weight $H$ of the  vertical edge in the boundary of the cell:
$$
\Sigma(H)= \{\eta\in\mathbb{C}: Re~\eta>0; 0<Im~\eta<H\}.
$$
This observation gives us the reconstruction rule for the curve $M$:
one has to glue $2\sharp\{vertical~edges\}$ half strips $\Sigma(H_s)$ in a way dictated by the 
graph. The conformal mapping of the glued surface to the sphere gives us a realization of the graph $\Gamma$, its vertexes with the odd value of $\ord(V)$ make up the branching set $\sf E$ of the curve $M$.  ~~~\bl

\subsection{Coordinate space of a graph}
The weights of an admissible graph $\{\Gamma\}$ have obvious linear restrictions. 
The dependencies arise from the symmetry of the graph as well as the conditions (W1), (W2) of Lemma \ref{five}.
Independent weights are given by the lengths $H(R)$ of the vertical edges $R$ in the closed upper half plane  and by the non-vanishing values of the width function $W(V)$ at the vertexes $V$ in closed upper half plane. 
They fill out an open convex polyhedron\footnote{Here and in what follows by the interior ($Int$) of a polyhedron we mean its
relative interior, that is the interior inside its affine hull.}   implicitly described in the last two items of Lemma \ref{five}. This polyhedron
is the product of the interior of the symplex $\Delta[\Gamma]$ spanned by positive variables $H(R)$, 
$~R\in Sk_1(\Gamma_\vert)\cap\hat{\mathbb{H}}$ with one normalization condition:
$$
\sum_{R\in Sk_1(\Gamma_\vert)\cap\mathbb{R}} H(R)+
2\sum_{R\in Sk_1(\Gamma_\vert)\cap\mathbb{H}} H(R) =\pi,
$$
and the interior of the cone $\mathfrak{C}[\Gamma]$ spanned by positive variables $W(V)$: $~V\in (Sk_0(\Gamma)\setminus Sk_0(\Gamma_\vert))\cap\hat{\mathbb{H}}$ that respect the partial order of nodes on the horizontal part of the graph given by the direction of arrows:
$$
  0<W(V_1)<W(V_2),
\qquad  if~ V_1<V_2.
$$

\begin{dfn}
The coordinate space of the graph $[\Gamma]$ we call the open polyhedron ${\cal A}[\Gamma]:=Int(\Delta[\Gamma]\times\mathfrak{C}[\Gamma])$. 
\end{dfn}

\begin{xmpl} 
For the graph $\Gamma$ depicted in the left panel of Fig. \ref{ExtGraph}, the genus $g[\Gamma]=2$, the number of ovals
$k[\Gamma]=1$. The coordinate space ${\cal A}[\Gamma]=\{H_1,H_2,H_3>0:~~2(H_1+H_2)+H_3=\pi\}\times$
$\{0<W_1<W_2\}$ has full dimension $dim {\cal A}[\Gamma]=dim {\cal H}_2^1=2g=4$.
\end{xmpl}

\begin{thrm}\cite{B2}, \cite{BBook}, {\rm \S 4.2.1 and \S 5.2~~~}
The mapping implied by Theorem \ref{CfromG} real analytically embeds the coordinate space  ${\cal A}[\Gamma]$ 
to the moduli space $\cal H$.
\end{thrm}

\begin{lmm}\cite{B2},\cite{BBook}, {\rm \S 4.2~~~}
The dimension of the cell  ${\cal A}[\Gamma]$ is not greater than $2g[\Gamma]$ and equals to  $2g[\Gamma]$ iff
the neighborhood of each vertex $V\in Sk_0(\Gamma)$ takes one of the following appearances

\begin{picture}(40,5)
\multiput(5,0)(15,0){3}{\circle*{1}}
\thicklines
\put(5,.5){\line(1,0){5}}
\put(5,-.5){\line(1,0){5}}
\put(11,0){,}

\put(15,.5){\line(1,0){10}}
\put(15,-.5){\line(1,0){10}}
\put(20,0){\vector(0,1){5}}
\put(26,0){,}

\put(30,0){\vector(1,0){5}}
\put(40,0){\vector(-1,0){5}}
\end{picture},

and vertices $V$ on the (real) symmetry axis additionally may be 
of the following vicinity types (up to central symmetry): \\

\begin{picture}(105,5)
\multiput(10,0)(15,0){6}{\circle*{1}}
\multiput(0,0)(15,0){7}{\dots}

\thicklines
\put(5,.5){\line(1,0){10}}
\put(5,-.5){\line(1,0){10}}
\put(10.5,5){\line(0,-1){10}}
\put(9.5,5){\line(0,-1){10}}

\put(25,.5){\line(1,0){5}}
\put(25,-.5){\line(1,0){5}}
\put(25,0){\vector(-1,0){5}}

\put(40,0){\vector(0,1){5}}
\put(40,0){\vector(0,-1){5}}
\thicklines
\put(35,.5){\line(1,0){10}}
\put(35,-.5){\line(1,0){10}}

\put(55.5,5){\line(0,-1){10}}
\put(54.5,5){\line(0,-1){10}}
\thinlines
\put(55,0){\vector(1,0){5}}
\put(55,0){\vector(-1,0){5}}

\put(70,0){\vector(-1,0){5}}
\put(70,0){\vector(1,0){5}}
\put(70,5){\vector(0,-1){5}}
\put(70,-5){\vector(0,1){5}}

\put(80,0){\dots}
\put(85,0){\vector(1,0){5}}
\put(85,5){\vector(0,-1){5}}
\put(85,-5){\vector(0,1){5}}

\put(95,-1){~, dotted line here is the real axis.}
\end{picture}
\label{FullDim}
\end{lmm}

\begin{rmk}
\label{NonSymm}
A graph $\Gamma$ may be associated to an arbitrary hypereliptic curve (not necessarily real) with two marked points interchanged by the involution $J$. The analogy of the above Lemma for such curves $M$ says that a graph $\Gamma(M)$ of a typical curve is composed only of the nodes with vicinities of the first three types of Lemma \ref{FullDim}. The presence of the mirror symmetry of the curve gives additional stable nodes on the symmetry axis. 
\end{rmk}

\begin{rmk}
A combinatorial algorithm listing all stable admissible graphs with given topological invariants $g,~k$ is available and will be published elsewhere. For the low dimensional components of the moduli space all codimension zero cells may be 
found manually just by trying all possible connections of the stable nodes listed above.
\end{rmk}

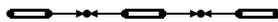
\begin{figure}[h!]
\centerline{
\begin{picture}(40,5)
\thicklines
\multiput(5,.5)(15,0){3}{\line(1,0){5}}
\multiput(5,-.5)(15,0){3}{\line(1,0){5}}
\thinlines
\multiput(5,0)(5,0){8}{\circle*{1}}
\multiput(10,0)(15,0){2}{\vector(1,0){5}}
\multiput(20,0)(15,0){2}{\vector(-1,0){5}}
\end{picture}}
\caption{\small  Graph $\Gamma$ encoding the unique cell of the moduli space ${\cal H}_2^3$.}
\label{H23}
\end{figure}

\begin{figure}[h!]
\begin{picture}(100,15)

\put(10,5){
\begin{picture}(40,20)
\thicklines
\put(5,.5){\line(1,0){10}}
\put(5,-.5){\line(1,0){10}}
\put(10.5,5){\line(0,-1){10}}
\put(9.5,5){\line(0,-1){10}}

\put(25,.5){\line(1,0){5}}
\put(25,-.5){\line(1,0){5}}
\thinlines

\multiput(5,0)(5,0){6}{\circle*{1}}
\put(10,5){\circle*{1}}
\put(10,-5){\circle*{1}}

\put(15,0){\vector(1,0){5}}
\put(25,0){\vector(-1,0){5}}
\end{picture}}

\put(50,5){
\begin{picture}(40,20)
\thicklines
\put(5.5,5){\line(0,-1){10}}
\put(4.5,5){\line(0,-1){10}}
\put(15,.5){\line(1,0){5}}
\put(15,-.5){\line(1,0){5}}

\put(30,.5){\line(1,0){5}}
\put(30,-.5){\line(1,0){5}}
\thinlines

\multiput(5,0)(5,0){7}{\circle*{1}}
\put(5,5){\circle*{1}}
\put(5,-5){\circle*{1}}

\put(5,0){\vector(1,0){5}}
\put(15,0){\vector(-1,0){5}}
\put(20,0){\vector(1,0){5}}
\put(30,0){\vector(-1,0){5}}
\end{picture}}

\put(90,5){
\begin{picture}(60,20)
\thicklines
\put(5,.5){\line(1,0){5}}
\put(5,-.5){\line(1,0){5}}
\put(20.5,5){\line(0,-1){10}}
\put(19.5,5){\line(0,-1){10}}

\put(30,.5){\line(1,0){5}}
\put(30,-.5){\line(1,0){5}}
\thinlines

\multiput(5,0)(5,0){7}{\circle*{1}}
\put(20,5){\circle*{1}}
\put(20,-5){\circle*{1}}

\put(10,0){\vector(1,0){5}}
\put(20,0){\vector(-1,0){5}}
\put(20,0){\vector(1,0){5}}
\put(30,0){\vector(-1,0){5}}
\end{picture}}
\end{picture}
\caption{\small  Graphs $[\Gamma_1]$, $[\Gamma_2]$, $[\Gamma_3]$ (left to right) and central symmetric graphs $[-\Gamma_1]$, $[-\Gamma_2]$ encode all full dimensional cells of the moduli space ${\cal H}_2^2$.}
\label{H22}
\end{figure}
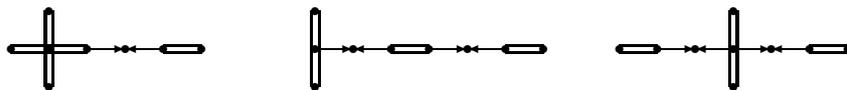

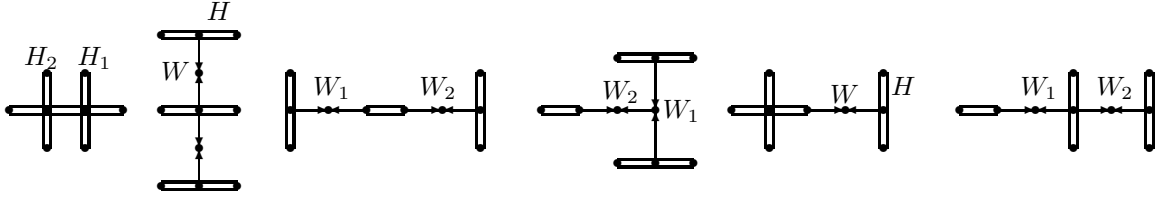
\begin{figure}[h!]
\begin{picture}(170,25)

\put(-5,15){
\begin{picture}(60,5)
\thicklines
\put(5,.5){\line(1,0){15}}
\put(5,-.5){\line(1,0){15}}
\put(10.5,5){\line(0,-1){10}}
\put(9.5,5){\line(0,-1){10}}
\put(15.5,5){\line(0,-1){10}}
\put(14.5,5){\line(0,-1){10}}
\thinlines

\multiput(5,0)(5,0){4}{\circle*{1}}
\put(10,5){\circle*{}}
\put(10,-5){\circle*{1}}
\put(15,5){\circle*{1}}
\put(15,-5){\circle*{1}}

\put(7,6){$H_2$}
\put(14,6){$H_1$}
\end{picture}}

\put(0,15){
\begin{picture}(60,10)
\thicklines
\put(20,.5){\line(1,0){10}}
\put(20,-.5){\line(1,0){10}}
\put(20,10.5){\line(1,0){10}}
\put(20,9.5){\line(1,0){10}}
\put(20,-10.5){\line(1,0){10}}
\put(20,-9.5){\line(1,0){10}}

\thinlines

\multiput(20,0)(5,0){3}{\circle*{1}}
\multiput(20,10)(5,0){3}{\circle*{1}}
\multiput(20,-10)(5,0){3}{\circle*{1}}
\put(25,5){\circle*{1}}
\put(25,-5){\circle*{1}}

\put(26,12){$H$}
\put(20,4){$W$}

\put(25,0){\vector(0,1){5}}
\put(25,10){\vector(0,-1){5}}
\put(25,0){\vector(0,-1){5}}
\put(25,-10){\vector(0,1){5}}
\end{picture}}

\put(17,15){
\begin{picture}(60,5)
\thicklines
\put(20.5,5){\line(0,-1){10}}
\put(19.5,5){\line(0,-1){10}}
\put(30,.5){\line(1,0){5}}
\put(30,-.5){\line(1,0){5}}
\put(45.5,5){\line(0,-1){10}}
\put(44.5,5){\line(0,-1){10}}
\thinlines

\multiput(20,0)(5,0){6}{\circle*{1}}
\put(20,5){\circle*{1}}
\put(20,-5){\circle*{1}}
\put(45,5){\circle*{1}}
\put(45,-5){\circle*{1}}

\put(20,0){\vector(1,0){5}}
\put(30,0){\vector(-1,0){5}}
\put(35,0){\vector(1,0){5}}
\put(45,0){\vector(-1,0){5}}

\put(23,2){$W_1$}
\put(37,2){$W_2$}

\end{picture}}

\put(50,15){
\begin{picture}(60,10)
\thicklines
\put(20,.5){\line(1,0){5}}
\put(20,-.5){\line(1,0){5}}
\put(30,7.5){\line(1,0){10}}
\put(30,6.5){\line(1,0){10}}
\put(30,-7.5){\line(1,0){10}}
\put(30,-6.5){\line(1,0){10}}

\thinlines

\multiput(20,0)(5,0){4}{\circle*{1}}
\multiput(30,7)(5,0){3}{\circle*{1}}
\multiput(30,-7)(5,0){3}{\circle*{1}}

\put(25,0){\vector(1,0){5}}
\put(35,0){\vector(-1,0){5}}

\put(35,7){\vector(0,-1){7}}
\put(35,-7){\vector(0,1){7}}

\put(36,-1){$W_1$}
\put(28,1.5){$W_2$}

\end{picture}}

\put(75,15){
\begin{picture}(60,5)

\thicklines
\put(25.5,5){\line(0,-1){10}}
\put(24.5,5){\line(0,-1){10}}

\put(40.5,5){\line(0,-1){10}}
\put(39.5,5){\line(0,-1){10}}
\put(20,.5){\line(1,0){10}}
\put(20,-.5){\line(1,0){10}}
\thinlines

\multiput(20,0)(5,0){5}{\circle*{1}}
\put(25,5){\circle*{1}}
\put(25,-5){\circle*{1}}
\put(40,5){\circle*{1}}
\put(40,-5){\circle*{1}}

\put(30,0){\vector(1,0){5}}
\put(40,0){\vector(-1,0){5}}

\put(41,2){$H$}
\put(33,1){$W$}

\end{picture}}

\put(105,15){
\begin{picture}(60,5)
\thicklines
\put(35.5,5){\line(0,-1){10}}
\put(34.5,5){\line(0,-1){10}}
\put(45.5,5){\line(0,-1){10}}
\put(44.5,5){\line(0,-1){10}}
\put(20,.5){\line(1,0){5}}
\put(20,-.5){\line(1,0){5}}
\thinlines

\multiput(20,0)(5,0){6}{\circle*{1}}
\put(35,5){\circle*{1}}
\put(35,-5){\circle*{1}}
\put(45,5){\circle*{1}}
\put(45,-5){\circle*{1}}

\put(25,0){\vector(1,0){5}}
\put(35,0){\vector(-1,0){5}}
\put(35,0){\vector(1,0){5}}
\put(45,0){\vector(-1,0){5}}

\put(28,2){$W_1$}
\put(38,2){$W_2$}

\end{picture}}
\end{picture}
\caption{\small  Graphs $[\Gamma_1], [\Gamma_2], \dots, [\Gamma_6]$ (left to right) and central symmetric graphs $[-\Gamma_4], [-\Gamma_5], [-\Gamma_6]$ encode all full dimensional cells of the moduli space ${\cal H}_2^1$.
The  weights of vertexes and vertical edges we denote as $W$ and $H$ respectively.}
\label{H21}
\end{figure}

\begin{xmpl}
The moduli spaces ${\cal H}_2^3$, ${\cal H}_2^2$ and ${\cal H}_2^1$ of dimension 4 will have respectively 1, 5 and 9
cells of full dimension; their graphs $[\Gamma]$ are shown in Figs. \ref{H23}, \ref{H22}, \ref{H21} up to central symmetry. 
The 6D spaces ${\cal H}_3^1$, ${\cal H}_3^2$, ${\cal H}_3^3$, ${\cal H}_3^4$ have respectively 24, 20, 7 and 1 codimension zero cells. 
Graphs encoding full dimensional cells of ${\cal H}_3^2$ are listed in \cite{BBook}, \S 4.2.
\end{xmpl}

\begin{rmk}
Several enumeration problems may be put forward: to find the number of  admissible graphs $[\Gamma]$ with given topological invariants $g,k$, and the number of full dimensional  admissible graphs (i.e. $\dim{\cal A}[\Gamma]=2g$). Same problems may be put forward for the curves without mirror symmetry.
\end{rmk}

\section{Polyhedral model of the Moduli space}
In the previous section we have built a decomposition of the total moduli space into smoothly embedded disjoint open polyhedra 
encoded by admissible (topological types of) graphs. We'll show that this decomposition is natural in the sense that
polyhedra of lower dimensions lie in the faces of higher dimensional polyhedra. This allows us 
to effectively build a PL model for each component of the total moduli space starting from the cells of codimension zero and gluing their faces according to certain identifications.

\subsection{Two types of faces of a coordinate space}
Each coordinate space ${\cal A}[\Gamma]$ is defined by the finite set of strict linear inequalities for the independent graph 
weights $H(R)$, $W(V)$ given in the conditions (W1)-(W2) of Lemma \ref{five}. On a face of the polyhedron certain inequalities 
turn into equalities. In other words, the weights of certain edges of the tree $\Gamma$ vanish and we want to interpret such 
degenerate weighted trees as regular trees $\Gamma'$, however with lower dimension of the coordinate space.

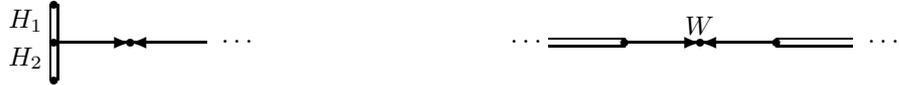
\begin{figure}[h!]
\begin{picture}(90,10)
\put(15,5){
\begin{picture}(20,10)
\thinlines
\put(4.5,5){\line(0,-1){10}}
\put(5.5,5){\line(0,-1){10}}
\thicklines

\multiput(5,5)(0,-5){3}{\circle*{1}}
\put(15,0){\circle*{1}}

\put(-1,2){$H_1$}
\put(-1,-3){$H_2$}
\put(5,0){\vector(1,0){10}}
\put(25,0){\vector(-1,0){10}}
\put(27,0){$\dots$}
\end{picture}}

\put(80,5){
\begin{picture}(20,10)
\thinlines
\put(5,-.5){\line(1,0){10}}
\put(5,.5){\line(1,0){10}}
\put(35,-.5){\line(1,0){10}}
\put(35,.5){\line(1,0){10}}

\thicklines
\multiput(15,0)(10,0){3}{\circle*{1}}
\put(0,0){$\dots$}
\put(47,0){$\dots$}
\put(23,1){$W$}
\put(15,0){\vector(1,0){10}}
\put(35,0){\vector(-1,0){10}}
\end{picture}}
\end{picture}

\caption{\small  Two typical junctions of branchpoints: $H_1+H_2\to0$ (left) and $W\to0$ (right). }
\label{BranchJunct}
\end{figure}

We have to distinguish between two types of faces in each polyhedron 
${\cal A}[\Gamma]$: 'outer' or 'exterior' faces correspond to the junction (= zero graph distance) of at least two branchpoints  -- see Fig. \ref{BranchJunct} -- and all the rest faces which we call 'inner' or 'interior'. 
\emph{Outer} faces  cannot be interpreted in terms of the same component of the moduli space\footnote{ Possibly, a more fruitful viewpoint is to consider points of outer faces as nodal curves living on the boundary of the moduli space component. Further desingularization allows to identify them with lower genus curves.} as shows  Proposition below. \emph{Inner} faces  may be further subdivided and identified with the lower dimension coordinate spaces of the same  moduli space component. Graphs $[\Gamma']$ corresponding to the inner faces of the coordinate space of graph $[\Gamma]$ we call subordinate to the latter graph and they may be obtained by successive application of two described below combinatorial procedures of \emph{Contraction} and \emph{Zipping} to the original graph.

\subsection{Elimination of zero weight edges: Definition}

Weighted graphs $\{\Gamma\}$ on the boundary of the coordinate space ${\cal A}[\Gamma]$ satisfy weak versions of axioms (W1)-(W2)
of Lemma \ref{five}:\\
{\it 
(W1*) Width function does not decrease along oriented edges and $W(V)=0$ if $V$ lies on the vertical part of the graph.\\
(W2*) The weights of vertical edges are nonnegative and their total sum is  $\pi$.}

\begin{rmk} Modified axioms (W1*)-(W2*) allow zero weight edges and surprisingly, may affect the topology of the graph $\Gamma$ too. For the graphs satisfying the original version of the axioms the value $d_{in}(V)$ is zero for vertexes $V$ of the vertical part of the graph. Indeed, starting in $V$ and moving against the arrows -- which is possible due to condition (T2) -- we eventually drop to other component of $\Gamma_\vert$. The value of the width function $W$ at this point should be strictly negative on the one hand and zero on the other hand due to the same axiom (W1). The weak version of the axiom admits
graphs with  positive $d_{in}(V)$ for $V\in Sk_0(\Gamma_\vert)$.
\end{rmk}

We shall apply two below described procedures \emph{Contraction} and \emph{Zipping} to weighted graphs corresponding to inner faces of a coordinate space, that is with positive graph distance between any two branchpoints.

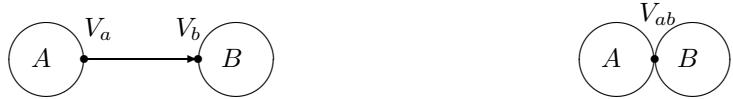
\begin{figure}[h]
\begin{picture}(40,10)
\put(25,0){\circle{10}}
\put(30,0){\circle*{1}}
\put(50,0){\circle{10}}
\put(45,0){\circle*{1}}
\put(100,0){\circle{10}}
\put(110,0){\circle{10}}
\put(105,0){\circle*{1}}
\put(30,3){$V_a$}
\put(42,3){$V_b$}
\put(23,-1){$A$}
\put(48,-1){$B$}
\put(103,5){$V_{ab}$}
\put(98,-1){$A$}
\put(108,-1){$B$}
\thinlines
\put(30,0){\vector(1,0){15}}
\end{picture}

\caption{\small  Left: the initial tree with zero weight horizontal edge $R=[V_a,V_b]$ and subtrees $A,B$ rooted in $V_a,V_b$ respectively;  Right: the edge $R$ is contracted.}
\label{Contr}
\end{figure}

{\bf Contraction} is applied to a zero weight horizontal edge $R$ of the graph $\Gamma$; its action on the graph is shown on Fig. \ref{Contr}. 

\begin{figure}
\begin{picture}(170,40)

\put(3,0){
\begin{picture}(100,40)
\thinlines
\put(20,10){\circle*{1}}
\put(20,20){\circle*{1}}

\put(20.2,10){\line(0,1){10}}
\put(19.8,10){\line(0,1){10}}
\put(5,13){$\Sigma$}

\put(20,10){\vector(-1,0){25}}
\put(15,10){\vector(-1,0){5}}

\put(25,10){\vector(1,0){5}}
\put(40,10){\vector(1,0){5}}
\put(55,10){\vector(1,0){5}}
\put(20,10){\vector(1,0){50}}

\put(5,10){\circle*{1}}

\put(35,10){\circle*{1}}
\put(50,10){\circle*{1}}
\put(65,10){\circle*{1}}

\put(15,11){$V_a$}
\put(33,11){$V_{a1}$}
\put(48,11){$V_{a2}$}
\put(63,11){$V_{a3}$}

\put(20,20){\vector(-1,0){25}}
\put(10,20){\vector(-1,0){5}}

\put(22,20){\vector(1,0){5}}
\put(40,20){\vector(1,0){5}}
\put(20,20){\vector(1,0){50}}

\put(0,20){\circle*{1}}

\put(33,20){\circle*{1}}
\put(50,20){\circle*{1}}

\put(15,16){$V_b$}
\put(31,16){$V_{b1}$}
\put(48,16){$V_{b2}$}

\put(5,5){\circle{10}}
\put(20,5){\circle{10}}
\put(35,5){\circle{10}}
\put(50,5){\circle{10}}
\put(65,5){\circle{10}}

\put(3,4){$A_{-1}$}
\put(18,4){$A$}
\put(33,4){$A_1$}
\put(48,4){$A_2$}
\put(63,4){$A_3$}

\put(0,25){\circle{10}}
\put(20,25){\circle{10}}
\put(33,25){\circle{10}}
\put(50,25){\circle{10}}
\put(-3,24){$B_{-1}$}
\put(18,24){$B$}
\put(31,24){$B_1$}
\put(48,24){$B_2$}

\put(15,11){$V_a$}
\put(15,16){$V_b$}
\end{picture}}
%
%
\put(82,5){
\begin{picture}(100,40)
\thinlines
\put(20,10){\circle*{1}}

\put(20,10){\vector(-1,0){25}}
\put(7,10){\vector(-1,0){5}}
\put(12,10){\vector(-1,0){5}}

\put(25,10){\vector(1,0){5}}
\put(32,10){\vector(1,0){5}}
\put(40,10){\vector(1,0){5}}
\put(55,10){\vector(1,0){5}}

\put(20,10){\vector(1,0){50}}

\put(0,10){\circle*{1}}
\put(5,10){\circle*{1}}
\put(33,10){\circle*{1}}
\put(38,10){\circle*{1}}
\put(50,10){\circle*{1}}
\put(65,10){\circle*{1}}

\put(5,5){\circle{10}}
\put(20,5){\circle{10}}
\put(38,5){\circle{10}}
\put(50,5){\circle{10}}
\put(65,5){\circle{10}}

\put(3,4){$A_{-1}$}
\put(18,4){$A$}
\put(36,4){$A_1$}
\put(48,4){$A_2$}
\put(63,4){$A_3$}

\put(0,15){\circle{10}}
\put(20,15){\circle{10}}
\put(33,15){\circle{10}}
\put(50,15){\circle{10}}
\put(-3,14){$B_{-1}$}
\put(18,14){$B$}
\put(31,14){$B_1$}
\put(48,14){$B_2$}
\end{picture}}
\end{picture}
\caption{\small  
Left: a zero height $H([V_a,V_b])=0$ strip with 
$0<W(V_{b1})<W(V_{a1})<W(V_{a2})=W(V_{b2})<W(V_{a3})$ is 'zipped' to a line (Right). 
Trees $A,~B$; $A_{\pm1},~B_{\pm1}$; $A2,~B2\dots$ are rooted at the outer sides of the strip.}
\label{Zipper}
\end{figure}
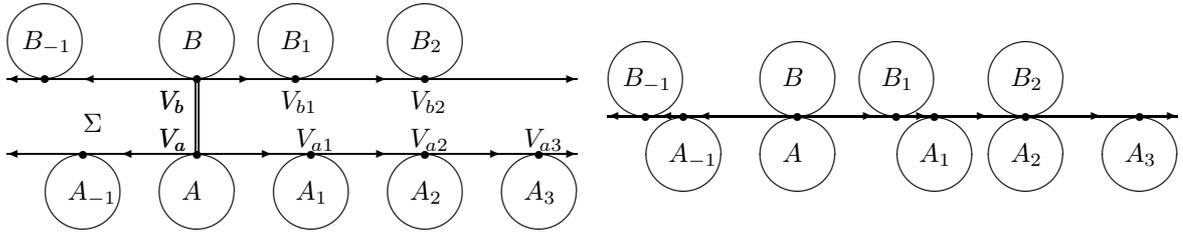

Elimination of a zero weight vertical edge of the graph is more complicated as it leads to the collapse of a whole strip in the decomposition of the plane considered in the proof sketch of Theorem \ref{CfromG} and therefore to deeper modification of the graph. {\bf Zipping} caused by a zero weight vertical edge $R$ of the graph is defined as follows. The edge is a vertical side of exactly two half-strips in the complement of the extended graph $Ext~\Gamma$. Each half-strip is collapsed to a ray so that points on the opposite sides with equal value of the width function $W$ are identified -- see Fig. \ref{Zipper}. Leaving for infinity horizontal edges of the modified extended graph  should be removed afterwards.

Informally speaking, the opposite sides of the strip passing through the zero weight vertical edge $R$ are considered as two pieces of a zipper fastener and attached one to the other in a natural way. The topology of the modified graph depends on the relations of the width function $W$ at the vertices on the opposite sides of the strip. It is exactly on this step that the faces of the coordinate space of the initial graph may be subdivided into smaller polyhedra. 

\begin{rmk}
As a byproduct of zipping, zero weight horizontal edges on the sides of the collapsed strip are contracted.
Their further contraction should be considered as a trivial action on the graph. 
\label{TrivCont}
\end{rmk}
\begin{rmk}
The procedure inverse to the elimination of zero weight edges leads to graphs with higher dimension of the coordinate space.
It consists in decay of unstable vertexes $V$ into more stable ones. The measure $codim(V)$ of a vertex $V$ stability  was introduced in \cite{B2}, \cite{BBook}, \S 4.1.3.
\end{rmk}

\subsection{Elimination of edges: Properties}

The application of the above two procedures leads -- in a finite number of steps --
to an admissible graph $\Gamma'$ as shows

\begin{lmm}
\label{WeakCons}
Contraction and zipping operations respect topological restrictions (T1)--(T3)
and weak restrictions on weights (W1*), (W2*) of the graph. 
The number of zero weight edges decreases after each operation.
\end{lmm}
{\bf Proof:}

(T1)
Both operations obviously preserve the dendritic nature of the graph: graphs in the right sides of the Figs. 
\ref{Contr}, \ref{Zipper} remain trees once $A,B,\dots$ are trees. Mirror symmetry of the weighted tree may be violated on a separate step, but eventually the symmetry will be repaired as two symmetric edges are collapsed or not simultaneously: see Lemma \ref{ElPerm} on permutability of operations.

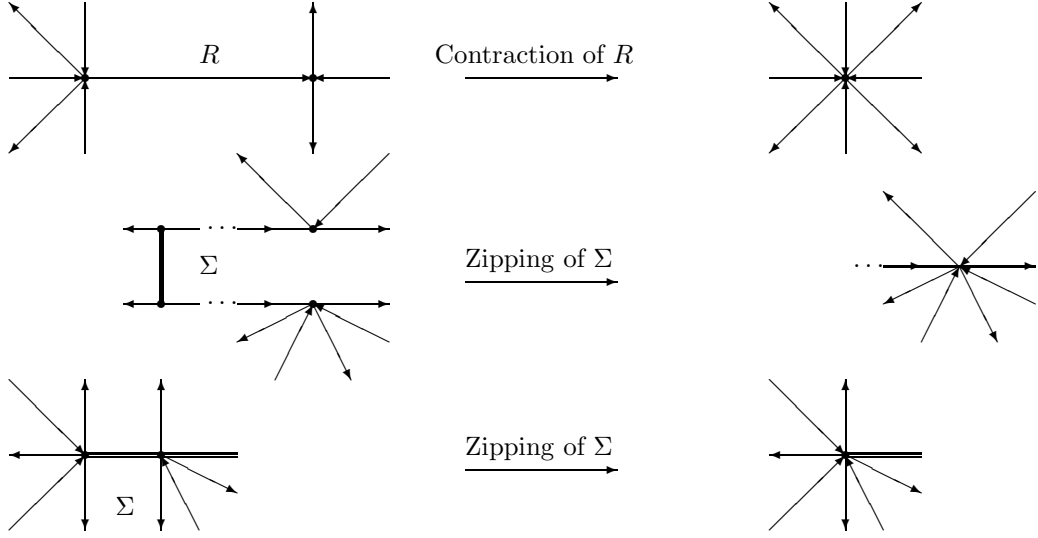
\begin{figure}
\begin{picture}(170,60)
\put(-5,50){  
\begin{picture}(70,20)
\put(15,10){\circle*{1}}
\put(45,10){\circle*{1}}

\put(15,10){\vector(1,0){30}}
\put(30,12){$R$}

\put(15,10){\vector(-1,-1){10}}
\put(15,10){\vector(-1,1){10}}
\put(15,20){\vector(0,-1){10}}
\put(5,10){\vector(1,0){10}}
\put(15,0){\vector(0,1){10}}

\put(45,10){\vector(0,-1){10}}
\put(45,10){\vector(0,1){10}}
\put(55,10){\vector(-1,0){10}}
\put(65,10){\vector(1,0){20}}
\put(61,12){Contraction of $R$}
\end{picture}}

\put(95,50){ 
\begin{picture}(70,20)
\put(15,10){\circle*{1}}

\put(15,10){\vector(-1,-1){10}}
\put(15,10){\vector(-1,1){10}}
\put(15,10){\vector(1,-1){10}}
\put(15,10){\vector(1,1){10}}

\put(25,10){\vector(-1,0){10}}
\put(5,10){\vector(1,0){10}}
\put(15,20){\vector(0,-1){10}}
\put(15,0){\vector(0,1){10}}

\end{picture}}

\put(15,20){ 
\begin{picture}(70,40)
\put(5,10){\circle*{1}}
\put(5,20){\circle*{1}}
\put(25,10){\circle*{1}}
\put(25,20){\circle*{1}}

\put(5.2,10){\line(0,1){10}}
\put(4.8,10){\line(0,1){10}}

\put(10,10){\vector(-1,0){10}}
\put(10,20){\vector(-1,0){10}}
\put(11,10){$\dots$}
\put(11,20){$\dots$}
\put(10,14){$\Sigma$}

\put(15,10){\vector(1,0){5}}
\put(15,20){\vector(1,0){5}}
\put(15,10){\vector(1,0){20}}
\put(15,20){\vector(1,0){20}}

\put(25,20){\vector(-1,1){10}}
\put(35,30){\vector(-1,-1){10}}

\put(25,10){\vector(1,-2){5}}
\put(25,10){\vector(-2,-1){10}}

\put(35,5){\vector(-2,1){10}}
\put(20,0){\vector(1,2){5}}

\put(45,15){Zipping of $\Sigma$}
\put(45,13){\vector(1,0){20}}
\end{picture}}

\put(100,25){ 
\begin{picture}(70,40)
\put(11,10){$\dots$}

\put(15,10){\vector(1,0){5}}
\put(15,10){\vector(1,0){20}}

\put(25,10){\vector(-1,1){10}}
\put(35,20){\vector(-1,-1){10}}

\put(25,10){\vector(1,-2){5}}
\put(25,10){\vector(-2,-1){10}}

\put(35,5){\vector(-2,1){10}}
\put(20,0){\vector(1,2){5}}
\end{picture}}

\put(0,0){ 
\begin{picture}(70,20)
\put(10,10){\circle*{1}}
\put(20,10){\circle*{1}}

\put(10,10.2){\line(1,0){20}}
\put(10,9.8){\line(1,0){20}}

\put(10,10){\vector(0,-1){10}}
\put(10,10){\vector(0,1){10}}
\put(10,10){\vector(-1,0){10}}
\put(0,0){\vector(1,1){10}}
\put(0,20){\vector(1,-1){10}}

\put(14,2){$\Sigma$}

\put(20,10){\vector(0,-1){10}}
\put(20,10){\vector(0,1){10}}
\put(20,10){\vector(2,-1){10}}
\put(25,0){\vector(-1,2){5}}

\put(60,10){Zipping of $\Sigma$}
\put(60,8){\vector(1,0){20}}
\end{picture}}

\put(100,0){
\begin{picture}(70,20)
\put(10,10){\circle*{1}}

\put(10,10.2){\line(1,0){10}}
\put(10,9.8){\line(1,0){10}}

\put(10,10){\vector(0,-1){10}}
\put(10,10){\vector(0,1){10}}
\put(10,10){\vector(-1,0){10}}
\put(0,0){\vector(1,1){10}}
\put(0,20){\vector(1,-1){10}}

\put(10,10){\vector(2,-1){10}}
\put(15,0){\vector(-1,2){5}}

\end{picture}}

\end{picture}
\caption{\small  Three typical junctions of nodes caused by contraction (on the top) and zipping (middle and bottom).}
\label{Junction}
\end{figure}

(T2) 
The horizontal edges leaving a vertex of the extended graph $Ext\Gamma$ alternate with the incident edges of other types -- incoming and vertical. Merging two vertexes during contraction of the edge or zipping conserves this property: see Fig. \ref{Junction}.  To return to the usual graph, we delete all edges leaving for infinity, therefore all the remaining 
outgoing edges are separated.

(T3) Let a contraction or zipping lead to the appearance of 
a vertex $V$ as a junction of several vertexes $V_a$, $V_b,\dots$. We show that 
\be
\label{AddOrd}
\ord(V)=\ord(V_a)+\ord(V_b)+\dots
\ee. 

For the contraction shown on Fig.\ref{Contr} we have a chain of equalities:
$$
\ord(V_{ab}):=
d_\vert^A(V_{ab})+d_\vert^B(V_{ab})+
2d_{in}^A(V_{ab})+2d_{in}^B(V_{ab})-2=
$$
$$
\bigg(d_\vert^A(V_a)+2d_{in}^A(V_a)-2\bigg)+
\bigg(d_\vert^B(V_b)+2(d_{in}^B(V_b)+1)-2\bigg)=
\ord(V_a)+\ord(V_b),
$$
where $d_\vert^A(V)$  stands for the degree of the vertex $V$ with respect to the vertical part of the subtree $A$,
$d_{in}^A(V)$ is the number of the incoming from the subtree $A$ edges  etc.

Zipping of a strip may lead to merging of more than two nodes. Preliminary contraction of zero weight edges on the horizontal sides of the strip
reduces the proof of the addition formula (\ref{AddOrd}) for the orders of vertexes to just two nodes junction. 

If $W(V_a)=W(V_b)=0$ as in left Fig. \ref{Zipper}, then 
$$
\ord(V_{ab}):=
d_\vert^A(V_{ab})+d_\vert^B(V_{ab})+
2d_{in}^A(V_{ab})+2d_{in}^B(V_{ab})-2=
$$
$$
\bigg((d_{\vert}(V_a)+1)+2d_{in}(V_a)-2\bigg)+
\bigg((d_{\vert}(V_b)+1)+2d_{in}(V_b)-2\bigg)=
\ord(V_a)+\ord(V_b).
$$

If otherwise $W(V_{a2})=W(V_{b2})>0$ as in left Fig. \ref{Zipper}, then 
$$
\ord(V_{a2b2}):=
2d_{in}^{A_2}(V_{a2b2})+2d_{in}^{B_2}(V_{a2b2})=
$$
$$
2(d_{in}^{A_2}(V_{a2})+1)-2 +2(d_{in}^{B_2}(V_{b2})+1)-2
=\ord(V_{a2})+\ord(V_{b2}).
$$

\begin{figure}

\begin{picture}(170,10)

\put(20,-5){ 
\begin{picture}(70,20)
\put(10,10){\circle*{1}}
\put(20,10){\circle*{1}}

\put(5,10.2){\line(1,0){20}}
\put(5,9.8){\line(1,0){20}}

\put(10,10){\vector(0,-1){5}}
\put(10,10){\vector(0,1){5}}
\put(20,10){\vector(0,-1){5}}
\put(20,10){\vector(0,1){5}}

\put(14,2){$\Sigma$}

\put(40,12){Zipping of $\Sigma$}
\put(40,10){\vector(1,0){20}}
\end{picture}}

\put(100,-5){
\begin{picture}(70,20)
\put(10,10){\circle*{1}}

\put(5,10.2){\line(1,0){10}}
\put(5,9.8){\line(1,0){10}}

\put(10,10){\vector(0,-1){5}}
\put(10,10){\vector(0,1){5}}
\end{picture}}
\end{picture}
\caption{\small  Zipping of a zero height strip leading to junction of two order zero vertexes.}
\label{Ord0}
\end{figure}
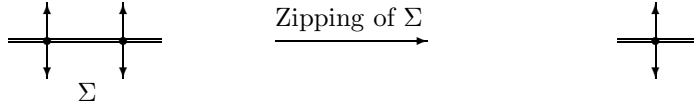

Suppose, the order of a newborn node is zero. We forbid junction of branchpoints, hence all orders in the above sum (\ref{AddOrd}) are even. 
Moreover, property (T2) implies that $\ord(V)>-2$, therefore  all the parent vertexes of the newborn
have zero order. According to (T3) property of the initial graph, each parent vertex is incident to two vertical and one or two outgoing 
edges only. This may only happen in case $Ext\Gamma$ fragment shown in Fig. \ref{Ord0} which after zipping gives a graph with (T3) property. 

Contracting and zipping change the topology of the graph, not weights. Therefore
the properties (W1*-W2*) remain intact. Contraction removes exactly one zero weight horizontal edge,
zipping eliminates one vertical edge of zero weight and all zero weight horizontal edges on the sides of the strip 
supported  by this vertical edge.  \bl

\begin{lmm}
Any two procedures eliminating zero weight edges -- contractions or zippings --  commute.
\label{ElPerm}
\end{lmm}
{\bf Proof.} Essentially, the fact follows from the local nature of both operations: 
nothing happens outside the edge we eliminate in case of contraction or outside the strip in case of zipping.
It follows from the (T1) property that two strips are either disjoint or have one edge of the extended graph in common.
See also Remark \ref{TrivCont}. ~~\bl

\begin{lmm}
 Elimination of zero weight edges does not merge the branchpoints.
\label{SameH}
\end{lmm}
{\bf Proof.} We apply the procedures of contraction and zipping only to graphs with positive graph distance 
between branchpoints.~~\bl  

The above three statements show that the outcome of the full chain of contractions/zippings applied to the 
weighted graph $\{\Gamma\}$ from the face of the coordinate space does not depend on the order of the operations 
and corresponds to an admissible weighted graph $\{\Gamma'\}$ from the same moduli space component. 
We call topological type of $\Gamma'$ subordinate to that of $\Gamma$ and write $[\Gamma']<[\Gamma]$.

\begin{lmm}
Once $[\Gamma']<[\Gamma]$, then the coordinate space ${\cal A}[\Gamma']$ is embedded
(possibly more than once, see Remark \ref{ToItself}) in the face of the polyhedron ${\cal A}[\Gamma]$.
\end{lmm}
{\bf Proof.} Subordination of the topological type of $\Gamma'$ to that of $\Gamma$ means that there is a chain of contractions/zippings
transforming the latter to the former. On each step the weights may be naturally and uniquely lifted to the 
larger graph. Eventually we get the weights on the graph $\Gamma$, but some of those will be zeros. This corresponds to a point on the boundary of the polyhedron $A[\Gamma]$. ~~\bl

\begin{dfn}
Dressing of a coordinate space ${\cal A}[\Gamma]$ we mean attaching all its interior faces:
$$
\hat{\cal A}[\Gamma]:=\mathop{\cup}\limits_{[\Gamma']\le[\Gamma]}{\cal A}[\Gamma'],
$$ 
the set is equipped with the ambient euclidean space topology. 
\end{dfn}

The basement for the described here incidence relations of the coordinate spaces of graphs is 

{\bf Proposition}
\label{Inc}
{\it The natural embedding of any dressed coordinate space $\hat{\cal A}[\Gamma]$ to the moduli space $\cal H$ is continuous and has no continuous extension to the exterior faces of ${\cal A}[\Gamma]$.}  

This statement is intuitively clear, however its rigorous proof requires special analytical techniques which is beyond the scope of this article. It will be given in a separate publication.  

\begin{figure}
\begin{picture}(170,20)
\unitlength=.8mm

\put(-10,0){
\begin{picture}(80,40)
\thinlines
\put(10,20.5){\line(1,0){15}}
\put(10,19.5){\line(1,0){15}}
\multiput(40,35.5)(0,-1){2}{\line(1,0){30}}
\multiput(40,5.5)(0,-1){2}{\line(1,0){30}}
\thicklines

\put(25,20){\vector(1,0){15}}
\put(55,20){\vector(-1,0){15}}
\put(55,35){\vector(0,-1){15}}
\put(55,5){\vector(0,1){15}}

\multiput(10,20)(15,0){4}{\circle*{2}}
\multiput(40,35)(15,0){3}{\circle*{2}}
\multiput(40,5)(15,0){3}{\circle*{2}}
\end{picture}}

\put(70,0){
\begin{picture}(80,40)
\thinlines
\put(10,20.5){\line(1,0){15}}
\put(10,19.5){\line(1,0){15}}
\multiput(25,35.5)(0,-1){2}{\line(1,0){30}}
\multiput(25,5.5)(0,-1){2}{\line(1,0){30}}
\thicklines
\put(25,20){\vector(1,0){15}}
\put(40,35){\vector(0,-1){15}}
\put(40,5){\vector(0,1){15}}

\multiput(10,20)(15,0){3}{\circle*{2}}
\multiput(25,35)(15,0){3}{\circle*{2}}
\multiput(25,5)(15,0){3}{\circle*{2}}
\end{picture}}

\put(155,0){
\begin{picture}(120,40)
\thinlines

\multiput(0,35)(0,-15){3}{
\multiput(0,0.5)(0,-1){2}{\line(1,0){15}}
\multiput(0,0)(15,0){2}{\circle*{2}}}

\thicklines
\put(15,20){\vector(1,0){15}}
\put(15,20){\oval(30,30)[r]}
\put(30,25){\vector(0,-1){5}}
\put(30,15){\vector(0,1){5}}
\put(30,20){\circle*{2}}
\end{picture}}
 
\end{picture}
\caption{\small  A  non-extendable chain of subordinate graphs $[\Gamma]>[\Gamma']>[\Gamma'']$ (left to right)
of the space ${\cal H}_2^1$ with the dimensions of their coordinate spaces $4,~3,~2$ respectively.}
\label{Chain}
\end{figure}
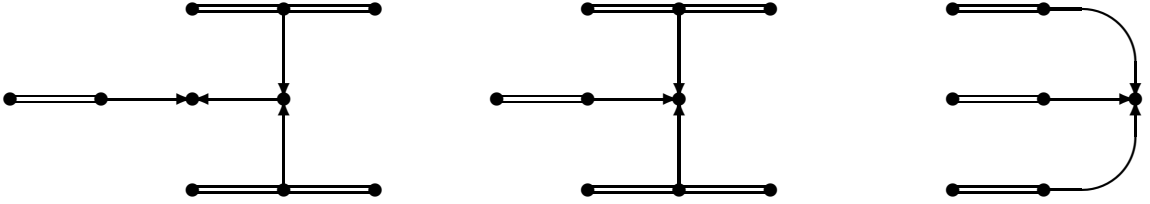

\subsection{Building Moduli and Labyrinth spaces}.

A piecewise-linear model of any component of the total moduli space ${\cal H}$ may be assembled in three steps:
\begin{itemize}
\item List all admissible graphs $[\Gamma]$ with full dimensional coordinate space ${\cal A}[\Gamma]$ (there are only finitely many of them once the genus is fixed). 
\item Dress each zero codimension coordinate space ${\cal A}[\Gamma]$.
\item Glue the dressed polyhedra $\hat{\cal A}[\Gamma]$ along their interior faces.
\end{itemize}
\begin{rmk}\ref{ToItself}
Note that on the last stage a full dimensional polyhedron may be glued to itself. Say if the graph $\Gamma$ contains the fragment shown on 
the left panel of Fig. \ref{BranchJunct} then two of its hyperfaces $\{H_1=0\}$ and $\{H_2=0\}$ correspond to the same subordinate graph $\Gamma'$ and should be glued one to the other.
\end{rmk}

\begin{figure}
\begin{picture}(170,30)
\unitlength=1.mm
\newsavebox{\Y}
\savebox{\Y}(40,25)[bl]{
\put(0,5){$\mathbb{H}$}
\multiput(0,0)(10,0){3}{\circle*{1}}
\multiput(0,20)(20,0){2}{\circle*{1}}
\put(10,10){\circle*{1}}
\multiput(20,0)(1,0){20}{.}

\put(0,.2){\line(1,0){20}}
\put(0,-.2){\line(1,0){20}}
\put(9.8,0){\line(0,1){10}}
\put(10.2,0){\line(0,1){10}}

\put(10.3,10){\line(1,1){10}}
\put(9.7,10){\line(1,1){10}}

\put(10.3,10){\line(-1,1){10}}
\put(9.7,10){\line(-1,1){10}}
}

\put(0,0){\usebox{\Y}} 
\qbezier(0,20)(0,25)(15,25)
\qbezier(15,25)(35,25)(35,0)
\qbezier(20,20)(30,20)(30,0)
\qbezier(10,10)(25,10)(25,0)

\put(50,0){\usebox{\Y}} 
\qbezier(70,20)(75,20)(75,0)

\qbezier(60,10)(60,25)(70,25)
\qbezier(70,25)(80,25)(80,0)

\qbezier(50,20)(50,30)(65,30)
\qbezier(65,30)(85,30)(85,0)

\put(100,0){\usebox{\Y}} 
\qbezier(120,20)(125,20)(125,0)

\qbezier(100,20)(100,25)(115,25)
\qbezier(115,25)(130,25)(130,0)

\qbezier(110,10)(95,10)(95,20)
\qbezier(95,20)(95,30)(115,30)
\qbezier(115,30)(135,30)(135,0)

\put(40,10){\vector(1,0){5}}
\put(40,12){$\beta_1$}

\put(90,10){\vector(1,0){5}}
\put(90,12){$\beta_2$}
\end{picture}
\caption{\small  Upper half of an exceptional graph $[\Gamma]$ with $g=3,~k=1$ admitting three canonical labyrinths transformed by standard generators $\beta_1$, $\beta_2$ of $Br_3$.}
\label{ExceptGraph}
\end{figure}
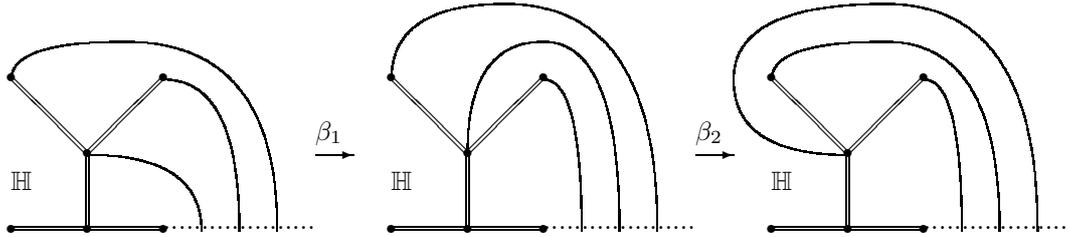

Gluing a labyrinth space from the dressed full dimensional polyhedra $\hat{\cal A}[\Gamma]:=\mathop{\cup}\limits_{\Gamma'\le\Gamma}{\cal A}[\Gamma']$ is somewhat more tricky. 
Lifting a point from the moduli space to the labyrinth space means attaching a labyrinth $\Lambda$
to the branching divisor $\sf E^+$ in the upper half plane.

\begin{dfn}
We call a labyrinth $\Lambda$ escorting the branching divisor $\sf E^+$ canonical iff 
it does not intersect the graph $\Gamma$ of the curve $M({\sf E})$ in the open upper half plane 
except at the endpoints of the arcs.
\end{dfn}

For \emph{non exceptional} graphs $[\Gamma]$ the canonical labyrinth is unique up to isotopy. Non exceptional here means that all branch 
points $V$ in the open upper half plane are the hanging vertexes of the graph (= stable: degree $d(V)=1$). In particular, all full 
dimensional graphs $[\Gamma]$ are non-exceptional. Exceptional graphs $[\Gamma]$ admit several\footnote{Namely $\prod_Vd(V)$, the product is 
taken over the branchpoints in the open upper half plane.}
canonical labyrinths $\Lambda$, $\Lambda'\dots$ which are interchanged by the action of braids: $\beta\cdot\Lambda:=\Lambda'$, 
$\beta\in Br$ -- see Fig. \ref{ExceptGraph}. The algorithm for the reconstruction of the braid from a couple of labyrinths escorting the same branchpoints may be found in \cite{DDRW}.

A coordinate space of a non-exceptional (e.g. full dimensional) graph $[\Gamma]$ has a canonical lift to the labyrinth space which we 
denote ${\cal LA}[\Gamma]$. The dressed polyherdon $\hat{\cal A}[\Gamma]=\cup_{\Gamma'\le\Gamma}{\cal A}[\Gamma']$ also has a canonical 
lift designated as ${\cal L}\hat{\cal A}[\Gamma]$, since all attached low dimensional faces inherit a  labyrinth   
from $\Gamma$ (which may be no longer canonical). All other lifts differ by the action of cover transformations of the universal covering 
represented by braids. Therefore, to build a PL model of the labyrinth space ${\cal L}_g^k$ we take all possible (finitely many)  
dressed full-dimensional polyhedra $\hat{\cal A}[\Gamma]$ corresponding to the given topological invariants $g,k$ and label 
them by braids $\beta\in Br_{g-k+1}$ (infinitely many if $g>k$).  Now we identify the common face  ${\cal A}[\Gamma_{12}]$ of 
two polyhedra  $\beta^1\cdot{\cal L}\hat{\cal A}[\Gamma_1]$ and $\beta^2\cdot{\cal L}\hat{\cal A}[\Gamma_2]$,  
$[\Gamma_{12}]<[\Gamma_1],[\Gamma_2]$, iff  $\beta^2=\beta^1\beta^{12}$, where the braid $\beta^{12}\in Br_{g-k+1}$ maps the  
labyrinth inherited by $\Gamma_{12}$ from $\Gamma_2$ to the labyrinth inherited from $\Gamma_1$. 

More formally, the total labyrinth space may be represented as follows:
\be
{\cal L}:=\mathop{\sqcup}\limits_{g,k}{\cal L}_g^k =
\mathop{\bigsqcup}\limits_{codim{\cal A}[\Gamma]=0, ~\beta\in Br}\beta\cdot{\cal L}\hat{\cal A}[\Gamma]/\sim,
\label{LRep}
\ee
with the following equivalence relation $~\sim~$ on the boundaries of the dressed coordinate spaces. 
If $[\Gamma']$ is subordinated to both $[\Gamma_1]$, $[\Gamma_2]$ then the coordinate space ${\cal A}[\Gamma']$ has two natural inclusions to the labyrinth space: as a face of ${\cal LA}[\Gamma_s]$, $s=1,2$. We identify ${\cal LA}[\Gamma'_1]$ with $\beta^{12}\cdot{\cal LA}[\Gamma'_2]$ where the braid $\beta^{12}$ maps the labyrinth of $\Gamma'$ inherited from $[\Gamma_2]$ to that inherited from $[\Gamma_1]$. The case $\Gamma_1=\Gamma_2$ here is not excluded: there may be more than one inclusion of the same coordinate space ${\cal A}[\Gamma']$ to a dressed coordinate space $\hat{\cal A}[\Gamma_1]$.

\section{PL model of Period map}\label{Sec5}

A point $M({\sf E})\in{\cal H}$ in the total moduli space  is presented by a (normalized) mirror-symmetric branching divisor $\sf E$.  Points of the covering labyrinth space above $M$ are distinguished by the choice of the labyrinth $\Lambda$ that escorts the set $\sf E^+$. A labyrinth determines a basis in the lattice of odd integer cycles on the (twice punctured at infinity) curve $M$ 
which is transported by Gauss-Manin connection -- see Lemma \ref{OddBasis}. The period mapping acting from each component of the labyrinth space to a suitable euclidean space is defined by formula (\ref{PM}) and it is evaluated at the distinguished basis of odd cycles on the surface. For a given point $\sf E^+\in {\cal H}$  we start with evaluation of the period map for a canonical labyrinth $\Lambda$ (i.e. not intersecting the graph of the curve $M(\sf E)$ in the upper half plane) which is unique for non-exceptional curves $M$.

\subsection{Period map restricted to a coordinate space}
To calculate the period map we need some preliminary considerations. 
One can single out a branch of the distinguished differential
$d\eta_M$ in any simply connected  domain of the plane avoiding the branch points of the curve $M$.
We take the complements to the labyrinth and to the graph in the open upper half plane as such sets.
In particular, we introduce a harmonic function 
$$
H'(x)=Im\int^x_*  d\eta_M,  
$$
in the complement  $\mathbb{H}\setminus\Lambda$ to the canonical labyrinth, with normalization $H'(x)=0$ for large real argument. The branch of the differential here is singled out by the condition: its residue at infinity equals $-1$. A harmonic function $H(x)$ conjugate to $W(x)$ is defined by the same formula  in the complement $\mathbb{H}\setminus\Gamma$ to the graph. The differentials $dH$ and $dH'$ coincide up to the sign in the components of the set $\mathbb{H}\setminus\{\Gamma\cup\Lambda\}$. This sign changes when we cross  either the labyrinth or the vertical edge of the graph $\Gamma$ and remains the same when we cross a horizontal edge. In particular, the sign $dH'/dH$ equals to $(-1)^{g+l}$ in the component of the complement $\mathbb{H}\setminus\{\Gamma\cup\Lambda\}$ bounded by the arcs $\Lambda_l$,  $\Lambda_{l+1}$ of the labyrinth, $l=k,\dots,g-1$. We call a component of the complement 'positive' or 'negative' when such is the sign of $dH'/dH$ in it. The boundary of the half plane cut along the graph $\Gamma$ has natural orientation and is divided into pieces labeled by vertical and horizontal edges of the graph. Cauchy-Riemann equations show that the boundary value $H(v)$ is locally constant at the 
pieces corresponding to horizontal edges and strictly decreases at the banks of the vertical edges. Therefore, the value of $H$
at a vertex $v$ of the boundary may be found as follows: we crawl along the boundary from $v$ toward $+\infty$ and sum up the weights $H(R)$ of all vertical edges we encounter (some of those may appear twice). 

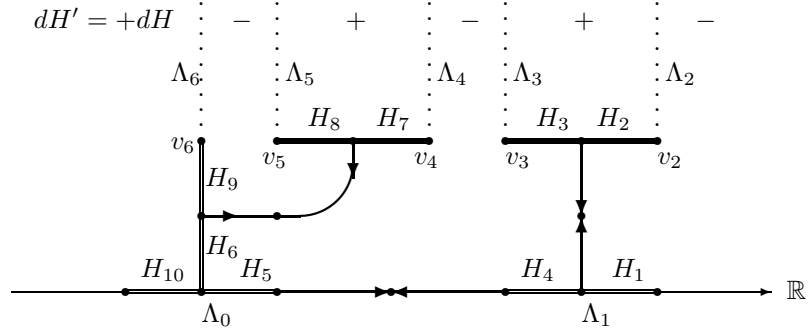
\begin{figure}
\hspace{3cm}
\begin{picture}(170,55)
\multiput(10,0)(10,0){3}{\circle*{1}}
\multiput(60,0)(10,0){3}{\circle*{1}}

\multiput(20,20)(10,0){7}{\circle*{1}}
\multiput(20,10)(10,0){2}{\circle*{1}}

\multiput(30,20)(20,0){2}{
\multiput(-0.5,0)(0,2){10}{.}}
\multiput(60,20)(20,0){2}{
\multiput(-0.5,0)(0,2){10}{.}}
\multiput(19.5,20)(0,2){10}{.}

\multiput(10,0.2)(0,-0.4){2}{\line(1,0){20}}
\multiput(60,0.2)(0,-0.4){2}{\line(1,0){20}}

\multiput(30,20.2)(0,-.4){2}{\line(1,0){20}}
\multiput(20.2,0)(-.4,0){2}{\line(0,1){20}}

\multiput(60,20.2)(0,-.4){2}{\line(1,0){20}}

\thicklines
\put(20,20){\oval(40,20)[br]}
\put(20,10){\vector(1,0){5}}
\put(40,20){\vector(0,-1){5}}

\put(30,0){\vector(1,0){15}}
\put(60,0){\vector(-1,0){15}}
\put(70,0){\vector(0,1){10}}
\put(70,20){\vector(0,-1){10}}

\put(70,10){\circle*{1}}
\put(45,0){\circle*{1}}

\thinlines
\put(80,0){\vector(1,0){15}}
\put(97,-1){$\mathbb{R}$}
\put(-5,0){\line(1,0){15}}
\put(20,-4){$\Lambda_0$}
\put(70,-4){$\Lambda_1$}

\put(16,28){$\Lambda_6$}
\put(31,28){$\Lambda_5$}
\put(51,28){$\Lambda_4$}
\put(61,28){$\Lambda_3$}
\put(81,28){$\Lambda_2$}

\put(16,19){$v_6$}
\put(28,17){$v_5$}
\put(48,17){$v_4$}
\put(60,17){$v_3$}
\put(80,17){$v_2$}

\put(-2,35){$dH'=+dH$}
\put(24,35){$-$}
\put(39,35){$+$}
\put(54,35){$-$}
\put(69,35){$+$}
\put(85,35){$-$}

\put(12,2){$H_{10}$}
\put(25,2){$H_5$}
\put(62,2){$H_4$}
\put(74,2){$H_1$}

\put(64,22){$H_3$}
\put(72,22){$H_2$}
\put(34,22){$H_8$}
\put(43,22){$H_7$}

\put(20.5,5){$H_6$}
\put(20.5,14){$H_9$}
\end{picture}
\caption{\small Upper half of a graph $\Gamma$, $M(\Gamma)\in{\cal H}_6^2$; canonical labyrinth $\Lambda=
(\Lambda_0,\dots,\Lambda_{6})$ (dotted lines in the upper half plane); weights of vertical edges $H_1,\dots,H_{10}$; 
auxiliary points $v_*$; real axis -- thin arrow}
\label{M62}
\end{figure}

\begin{xmpl}
Boundary values of the function $H(v)$ for curve $M(\Gamma)$ with the graph shown in Fig. \ref{M62} are as follows 
$$
\begin{array}{l}
H(v_2)=H_1+H_2;\qquad
H(v_6)=H_1+2H_2+2H_3+H_4+H_5+H_6+2H_7+2H_8+H_9;\\
H(v_3)=H_1+2H_2+H_3;\qquad
H(v_5)=H_1+2H_2+2H_3+H_4+H_5+H_6+2H_7+H_8;\\ 
H(v_4)=H_1+2H_2+2H_3+H_4+H_5+H_6+H_7.
\end{array}
$$
\end{xmpl}

\begin{lmm}\cite{B2}, \cite{BBook}, {\rm \S 5.3~~}
Periods map $\Pi_s:=\langle\Pi|C_s\rangle:=i\int_{C_s}d\eta_M$ evaluated at the odd cycles $C_s$
determined by the canonical labyrinth $\Lambda$ is given by the expression:
\be
\Pi_s(\Lambda)= \left\{
\begin{array}{ll}
2\sum\limits_{R\subset\Lambda_s} \epsilon(R)H(R),&\qquad s=0,\dots,k-1;\\
4(-1)^{s+g}H(v_s),&\qquad s=k,\dots,g.
\end{array}
\right.
\label{Ps}
\ee
where the sum in the formula for $s=0,\dots,k-1$ is taken over all real vertical edges $R$ of the graph
$\Gamma$ which make up the segment $\Lambda_s$ and $\epsilon(R)$ is the sign of the component
of $\mathbb{H}\setminus(\Lambda\cup\Gamma)$ bordered by $R$; $v_s$ in the other formula for $s=k,\dots,g$, is the meeting points 
of the graph and the $s-$th arc of the labyrinth in the upper half plane, it lies on the boundary of $\mathbb{H}\setminus\Gamma$.
\label{Permap}
\end{lmm}

{\bf Proof.} Each of the basic cycles $C_s$ changes its orientation after complex conjugation, hence we can integrate along upper half of that cycle:
$$
\Pi_s:=-Im~\int_{C_s}d\eta_M=-2Im~\int_{C_s\cap\mathbb{H}}d\eta_M=-2\int_{C_s\cap\mathbb{H}}dH'.
$$
We consider  separately two cases: 
(a) $s=0,\dots,k-1$ when the integration path is the upper bank of the real segment
$\Lambda_s$ and 
(b) $s=k,\dots,g$ when we integrate along both banks of the cut $\Lambda_s$ drawn in the upper half plane.

(a) We integrate $dH'$ along the upper bank of the real segment $\Lambda_s$
$$
\Pi_s/2=-\int_{\Lambda_s} dH'=
-\sum_{R\subset\Lambda_s}\int_R dH'=
-\sum_{R\subset\Lambda_s}\pm\int_R dH=  
\sum_{R\subset\Lambda_s}\epsilon(R)H(R),  
$$
where $\epsilon(R)$ in the last formula is the sign of the component of the complement 
to the graph and the labyrinth whose boundary contains the edge $R$.

(b) We integrate $dH'$ along both banks of the cut $\Lambda_s$. Instead, we can integrate along the left bank  and double the result. This integration path lies in the component of the complement to both labyrinth and the graph $\Gamma$  bounded by the arcs $\Lambda_s$ and $\Lambda_{s-1}$. We integrate $(-1)^{s+g+1}dH$ instead of $dH'$:
$$
(-1)^{s+g}\Pi_s/4=H(v_s)=\sum_{R>v_s} H(R)
$$
where $v_s$ is the point of touching of the arc $\Lambda_s$ and the graph; the sum is taken over all vertical edges of the boundary of the upper half-plane cut along the graph that lie  between $v_s$ and $+\infty$.
~~\bl

\begin{rmk}
We see that the period map depends only on the coordinates $H$ of the simplicial factor $\Delta[\Gamma]$ of the space ${\cal A}[\Gamma]$. 
If it does not lead to a conclusion (e.g. in case of non-exceptional $[\Gamma]$) we consider the period mapping (\ref{Ps}) directly from the 
simplex $\Delta[\Gamma]$ to the euclidean space and designate it by the same letter $\Pi$.
\end{rmk}

\begin{xmpl}
The period map for curve $M\{\Gamma\}$ with the graph $[\Gamma]$ shown in Fig. \ref{M62} equals to
$$
\begin{array}{l}
\Pi_0=2(H_{10}-H_5);\qquad \Pi_6=4(H_1+2H_2+2H_3+H_4+H_5+H_6+2H_7+2H_8+H_9);\\
\Pi_1=-2(H_4+H_1);\qquad \Pi_5=-4(H_1+2H_2+2H_3+H_4+H_5+H_6+2H_7+H_8);\\ 
\Pi_2=4(H_1+H_2);\qquad \Pi_4=4(H_1+2H_2+2H_3+H_4+H_5+H_6+H_7);\\
\Pi_3=-4(H_1+2H_2+H_3).\\
\end{array}
$$
One can check that $\sum_{s=0}^6\Pi_s=2(H_1+H_4+H_5+H_{10})+4(H_2+H_3+H_6+H_7+H_8+H_9)=2\pi$
due to the normalization of weights of vertical edges.
\end{xmpl}

The equality $\sum_{s=0}^g\Pi_s=2\pi$  follows from the fact that the sum of all basic odd cycles is homological to the small circle encompassing the pole of the distinguished differential $d\eta_M$. This means that the image of the period map lies in the hyperplane of $\mathbb{R}^{g+1}$. Let us check the equality by a straightforward calculation:
$$
\sum\limits_{s=k}^g\Pi_s=4(H(v_g)-H(v_{g-1}))+4(H(v_{g-2})-H(v_{g-3}))+\dots (-1)^{g+k}4H(v_k)
=4\sum_R H(R),
$$
the last sum is taken over the vertical edges from the boundaries of all 'negative' domains in decomposition of $\mathbb{H}$
by the graph and the labyrinth.  
Adding $\sum_{s=0}^{k-1}\Pi_s$ to the latter sum, we get a quadruple sum of the weights $H(R)$ of all 
vertical edges in the open upper half plane and the double sum of the weights of the graph  vertical edges  in the
real axis. The total value equals $2\pi$ due to the weight normalization condition (W2).

\subsection{Image of coordinate space and local fibers of period map} 

Lemma \ref{Permap} allows us to calculate the range of the periods map.
We know that any labyrinth space is tiled by the translations of the copies of coordinate spaces ${\cal A}[\Gamma]$.
A coordinate space is embedded into the component of the moduli space and then lifted to its covering space by drawing  
a (usually unique) canonical labyrinth. Lemma \ref{Permap} computes the period mapping on the latter cell; the mapping in the  
translated cells is related to the computed one by the Burau representation (\ref{BuRep}) -- see Theorem \ref{PMfacts}.
We know that any coordinate space is the product of the interior of the simplex spanned by the \emph{vertical} variables $H(R)$ and the interior of the cone spanned by the \emph{horizontal} variables $W(V)$. The period map depends on the weights of the vertical edges only, so the image of a coordinate space under the period mapping may be obtained as follows. Evaluate the period mapping at the vertexes of the simplicial factor of the coordinate space
and take the relative interior of the convex hull of thus obtained points. 

Following this procedure, the range of the periods map was calculated for all values of topological invariants $g,k$ in \cite{B2}, \cite{BBook}, \S.5.3. Here we reproduce the result for the genus $g=2$ and $k=1,2,3$, see Fig. \ref{ImGen2}.
For the purposes of the next section we also calculate the fibers of the period map inside each full dimensional cell of the
moduli space corresponding to genus $g=2$.

\subsubsection{Three real ovals} \label{FiberH23}
The image of the only coordinate space from ${\cal L}_2^3={\cal H}_2^3$  projected to the plane of variables $\Pi_0,\Pi_1$
is the open triangle:
\be
\label{PIH23}
\Pi_0>0;\quad\Pi_1>0;\quad\Pi_0+\Pi_1<2\pi.
\ee
The fibers of this map are quadrants $\mathbb{R}^+_2$ spanned by the weights $W_1$, $W_2$ of two vertexes 
from the horizontal subgraph -- see Fig \ref{H23}. 

\subsubsection{Two real ovals} 
The space ${\cal L}_2^2={\cal H}_2^2$ contains five full dimensional polyhedra corresponding to the topological graphs
$[\Gamma_1]$, $[\Gamma_2]$, $[\Gamma_3]$ of Fig. \ref{H22} (left to right) and also their central symmetric graphs  
$[-\Gamma_1]$, $[-\Gamma_2]$. Simplicial factor $\Delta[\Gamma]$ of the coordinate space ${\cal A}[\Gamma]$ is 3-dimensional 
for $[\Gamma]=[\pm\Gamma_1]$ and 2-dimensional in the remaining three cases. Hence the fiber of the period map  restricted to 
coordinate spaces ${\cal A}[\pm\Gamma_1]$ will be a half-strip: the product of an interval (section of the tetrahedron by a line) 
times a ray (conical factor of the coordinate space). The fibers in three other coordinate spaces corresponding to 
$[\Gamma]=[\pm\Gamma_2],[\Gamma_3]$, will be their conical factors: open quadrants spanned by two positive weights of vertexes in the 
horizontal subgraphs. Images of full dimensional coordinate spaces projected to the plane of coordinates $\Pi_0, \Pi_1$ as well as the 
local fibers of the period map in those polyhedra is given in the Table \ref{Tab22} below.

\begin{table}[h]
\begin{tabular}{l|l|l}
 $\Gamma$ & Image $\Pi({\cal LA}[\Gamma])$ & Fiber $\Pi^{-1}(\Pi^*$), $\Pi^*\in$ Image\\
\hline
$\Gamma_1$& $Int(a_+\cup b)$& $(H,W):~~0<4H<2\pi-|\Pi_0^*|+\Pi_1^*;0<W$\\
$-\Gamma_1$& $Int(a_-\cup b)$& $(H,W):~~0<4H<2\pi-\Pi_0^*-|\Pi_1^*|;0<W$\\
$\pm\Gamma_2$& $Int(a_\pm)$& $(W_1,W_2):~0<W_1,W_2$\\
$\Gamma_3$& $Int(b)$& $(W_1,W_2):~0<W_1,W_2$
\end{tabular}
\caption{\small  Images and fibers of the periods map in full dimensional cells of ${\cal H}_2^2$. 
Here $a_\pm ,b$ are three closed triangles shown in the Fig. \ref{PIH22}; $H$ is the weight of the only vertical edge in the upper half plane;  $W_*$ are the weights of the vertexes in the horizontal subgraph.}
\label{Tab22}
\end{table}

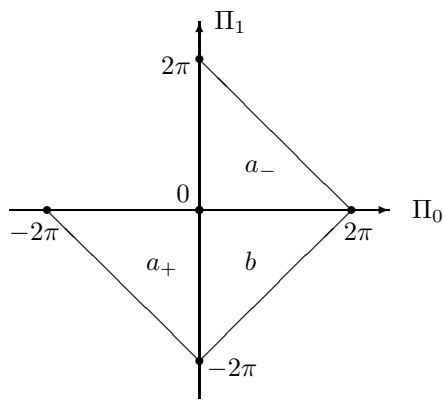
\begin{figure}[h!]
\hspace{5cm}
\begin{picture}(50,60)
\unitlength=1mm
\multiput(25,5)(0,20){3}{\circle*{1}}
\multiput(5,25)(20,0){3}{\circle*{1}}

\put(25,0){\vector(0,1){50}}
\put(27,49){$\Pi_1$}
\put(26,3){$-2\pi$}
\put(20,43){$2\pi$}
\put(22,26){$0$}

\put(0,25){\vector(1,0){50}}
\put(53,24){$\Pi_0$}
\put(0,21){$-2\pi$}
\put(44,21){$2\pi$}

\put(45,25){\line(-1,1){20}}
\put(25,5){\line(1,1){20}}
\put(25,5){\line(-1,1){20}}

\put(18,17){$a_+$}
\put(31,17){$b$}
\put(31,30){$a_-$}
\end{picture}
\caption{\small  The range of the period mapping in the coordinate spaces of ${\cal H}_2^2$.
Black dots are the images of the vertexes of the simplexes $\Delta[\Gamma]$} 
\label{PIH22}
\end{figure}

\subsubsection{One real oval} 
The moduli space ${\cal H}_2^1$ contains nine full dimensional polyhedra corresponding to the topological graphs
$[\Gamma_1]$, $\dots$, $[\Gamma_6]$ of Fig. \ref{H21} (left to right) and also their central symmetric graphs  
$[-\Gamma_4]$, $[-\Gamma_5]$, $[-\Gamma_6]$. The dimension of the simplicial factor $\Delta[\Gamma]$ of the coordinate space 
${\cal A}[\Gamma]$ is four when $[\Gamma]=[\Gamma_1]$, three for $[\Gamma]=[\Gamma_2],[\pm\Gamma_5]$ and two in the remaining three cases. 
The fiber of the period map  restricted to a  coordinate space ${\cal LA}[\Gamma_1]$ -- a section of open 4-simplex by a 2-plane -- 
may be either a pentagon, a quadrilateral, or a triangle. Period map fiber inside coordinate space for $[\Gamma]=[\Gamma_2],[\pm\Gamma_5]$ 
is a half-strip, the product of an interval (section of tetrahedron by a line) and a ray (conical factor of coordinate space). 
The fibers in five other coordinate spaces corresponding to  $[\Gamma]=[\Gamma_3],~[\pm\Gamma_4],~[\pm\Gamma_6]$, will be their conical 
factors spanned by two positive weights of vertexes in the horizontal subgraphs. Images of full dimensional coordinate spaces projected to 
the plane of coordinates $\Pi_1, \Pi_2$ as well as the fibers of the period map in those polyhedra are given in the Tab. \ref{Tab21}.

\begin{table}[h!]
\begin{tabular}{l|l|l}
 $\Gamma$ & Image $\Pi({\cal LA}[\Gamma])$ & Fiber $\Pi^{-1}(\Pi^*$), $\Pi^*\in$ Image\\
\hline
$~\Gamma_1$& $Int(a\cup b\cup c_+\cup c_-\cup d)$& Rectangle, $\Pi^*\in Int(a)$\\
&& Pentagon, $\Pi^*\in Int(b)$\\
&& Trapezoid, $\Pi^*\in Int(c_+\cup c_-)$\\
&& Triangle $\Pi^*\in Int(d)$\\
$~\Gamma_2$& $Int(b\cup c_+\cup c_-\cup d)$& Half -strip $(H,W):~~min(-\Pi^*_1,\Pi_1^*+\Pi_2^*)>4H$\\
&& $>max(0,2\Pi_2^*+\Pi_1^*-4\pi);\quad W>0$\\
$~\Gamma_3$& $Int(a)$& Quadrant $(W_1,W_2)\in\mathbb{R}_+^2$\\
$\pm\Gamma_4$& $Int(b\cup c_\pm)$& Cone $(W_1,W_2): 0<W_1<W_2$\\
$\pm\Gamma_5$& $Int(a\cup b\cup c_\pm)$& Half-strip\\
$\pm\Gamma_6$& $Int(b\cup c_\pm)$&  Quadrant $(W_1,W_2)\in\mathbb{R}_+^2$\\
\end{tabular}
\caption{\small  Images and fibers in ${\cal H}_2^1$.
Here $a,b,c_\pm,d$ are closed triangles shown in the left panel of Fig. \ref{PIH21}; $H$ is the weight of the vertical edge in the upper half plane;  $W_*$ are the weights of the vertexes in the horizontal subgraph -- see Fig \ref{H21}.
}
\label{Tab21}
\end{table}

\begin{figure}
\hspace{.5cm}
\begin{picture}(50,55)
\unitlength=1mm
\multiput(5,45)(20,0){3}{\circle*{1}}
\put(45,5){\circle*{1}}
\multiput(45,5)(0,20){3}{\circle*{1}}
\multiput(42,16)(-8,8){4}{\circle{1}}

\put(45,0){\vector(0,1){50}}
\put(40,48){$\Pi_2$}
\put(46,44){$4\pi$}
\put(46,24){$2\pi$}
\put(46,2){$0$}

\put(0,5){\vector(1,0){53}}
\put(56,4){$\Pi_1$}
\put(1,0){$-4\pi$}
\put(21,0){$-2\pi$}

\put(45,5){\line(-1,1){40}}
\put(45,5){\line(-1,2){20}}
\put(45,25){\line(-1,1){20}}
\put(45,25){\line(-2,1){40}}
\put(45,45){\line(-1,0){40}}

\put(5,4){\line(0,1){2}}
\put(25,4){\line(0,1){2}}

\put(38,38){$a$}
\put(33,33){$b$}
\put(27,27){$d$}
\put(40,20){$c_+$}
\put(20,40){$c_-$}
\end{picture}
\hspace{2cm}
\begin {picture}(70,50)
\unitlength=1mm
\put(10,5){\circle*{1}}

\put(5,5){\vector(1,0){60}}
\put(10,35){\line(1,0){37}}
\put(45,5){\line(0,1){32}}
\put(55,5){\line(-1,1){35}}
\put(55,5){\circle*{1}}
\put(45,5){\circle*{1}}
\put(45,15){\circle*{1}}
\put(25,35){\circle*{1}}
\put(10,35){\circle*{1}}

\put(67,3){$4H_1$}
\put(11,6){$0$}
\put(51,0){$\Pi_1^*+\Pi_2^*$}
\put(40,0){$-\Pi_1^*$}

\put(10,0){\vector(0,1){40}}
\put(5,42){$4H_2$}
\put(-5,33){$4\pi-\Pi^*_2$}
\end {picture}
\caption{\small  Left: range of the periods map in (canonically lifted) stable coordinate spaces of ${\cal H}_2^1$.
Black dots are the images of extreme points of the simplexes $\Delta[\Gamma]$; white dots -- the orbit of $Br_2$ inside the large triangle $\Delta$. ~~~~~~ Right: fiber of the periods map restricted to ${\cal LA}[\Gamma_1]$.} 
\label{PIH21}
\end{figure}
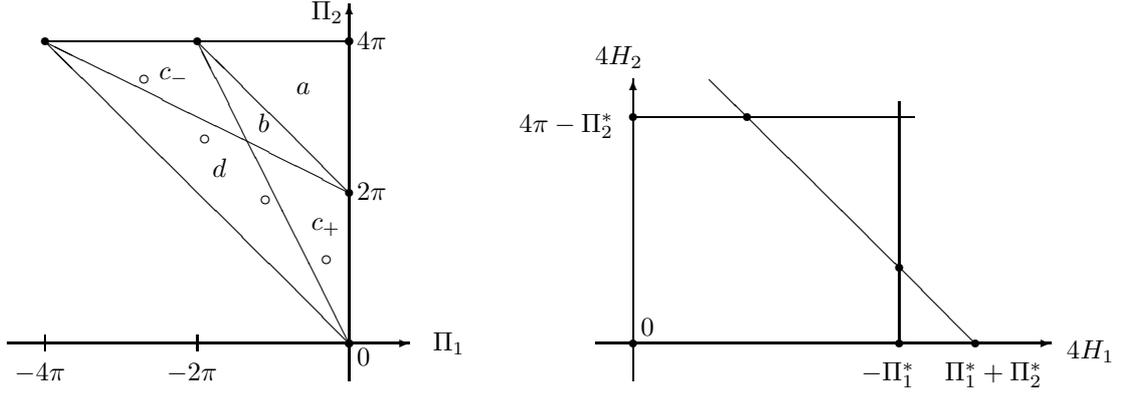

\begin{xmpl}
Let us consider the section of the 4-simplex ${\cal A}[\Gamma_1]$ in greater detail.
The fiber $\Pi^{-1}(\Pi^*)$ of the period map we parametrize by the weights $H_1,H_2$ of two 
vertical edges of the graph in the upper half plane (see most left picture in Fig. \ref{H21}). 
The positivity of weights of three vertical edges of $[\Gamma_1]$ on the real (symmetry) axis imply
the following inequalities:
$$
\begin{array}{lll}
 0<&4H_1&<-\Pi_1^*,\\
 &4(H_1+H_2)&<\Pi_1^*+\Pi_2^*,\\
 0<&4H_2&<4\pi-\Pi_2^*.\\
\end{array}
$$
Depending on the relations between components of $\Pi^*$ we get either a pentagon or a rectangle or a trapezoid or a triangle
-- see the right panel of Fig. \ref{PIH21}.
\end{xmpl}

The braid group $Br_2=\mathbb{Z}$ acts as cover transformations on the labyrinth space ${\cal L}_2^1$.
The images of non-canonical lifts of the coordinate spaces are related to the images we just calculated by the action of 
Burau representation which is generated by the matrix 
$$
B:=
\left(
\begin{array}{c}
~0-1\\1~~~2
\end{array}
\right)
$$

\begin{figure}
\vspace{-1mm}
\centerline{\includegraphics[scale=.9]{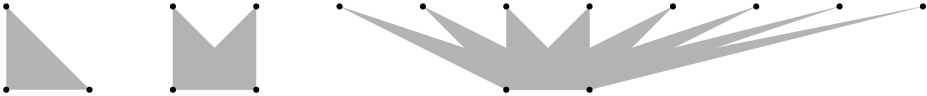}}
\vspace{-1mm}
\caption{\small  Images $\Pi({\cal L}_2^k)$ with  $k=3,2,1$ (left to right and only partially for $k=1$).}
\label{ImGen2}
\end{figure}

Part of the full image of the labyrinth space ${\cal L}_2^1$ is shown on the most right picture in Fig. \ref{ImGen2}.
It is the union of images of the large triangle $\Delta=a\cup b\cup c_\pm\cup d$ shown in the left panel of Fig. \ref{PIH21}
under the action of integer powers of the matrix $B$.

\subsection{Fibers on the boundary of coordinate space} 
Typically, a fiber of the periods map transversally intersects the boundaries of the full dimensional cells which tile a labyrinth space. 
However, some exceptional fibers may (locally) lie in the boundary. Here we specify the conditions of this relatively rare phenomenon.

The period mapping inside the dressed coordinate space $\hat{\cal A}[\Gamma]$ lifted to the labyrinth space depends only on the coordinates $H(\cdot)$ in the simplex $\Delta[\Gamma]$. The conical factor $\mathfrak{C}[\Gamma]$ of the coordinate space will always be a factor of the local fiber. Hence, we have to study the conditions when a linear map $\Pi$ restricted to a simplex $\Delta$ has the fiber $f$ of the maximal dimension lying in the boundary of the simplex.

The fiber $f$ belongs to the translated kernel of the linear map $\Pi$, hence the inequality 
\be
\dim\Delta-\dim\Pi\Delta\ge \dim f,
\label{DimEq}
\ee
which becomes the equality for instance when the fiber $f$ intersects the (relative) interior of the simplex.

\begin{lmm}
\label{fDelta}
The fiber $f$ has (maximal) dimension equal to the left hand side of (\ref{DimEq})
iff there exists a boundary simplex $\Delta'\le\Delta$ such that the fiber $f$ intersects the interior\footnote{
By the interior of 0-simplex we mean the point itself?}  of $\Delta'$
and the codimension of this simplex (in $\Delta$) equals to the codimension of it's image $\Pi\Delta'$ (in $\Pi\Delta$) 
\end{lmm}
{\bf Proof} 1. Suppose we've found the boundary simplex $\Delta'\le\Delta$ with the required property:
$$
\dim\Delta - \dim\Delta' =\dim\Pi\Delta - \dim\Pi\Delta',
$$
with the fiber $f$ intersecting the relative interior of $\Delta'$. Consider $f':=f\cap\Delta'$ -- the fiber of $\Pi$
in smaller simplex. The inequality (\ref{DimEq}) now becomes 
\be
\dim\Delta'-\dim\Pi\Delta'=\dim f',
\label{DimEq2}
\ee
Two latter equalities together with obvious the inequality $\dim f'\le\dim f$ give the inequality opposite to (\ref{DimEq}).
 
2. Suppose now we have an equality sign in (\ref{DimEq}). Let $\Delta'\le\Delta$
be the largest simplex in the boundary such that $f$ intersects its interior (take the convex hull of all simplexes with this property). This simplex contains the whole fiber $f$. Indeed, a point from $f\setminus \Delta'$ belongs to the interior of a unique simplex $\Delta''\le\Delta$
and the convex hull of $\Delta'$ and $\Delta''$ is larger than $\Delta'$ or $\Delta''\le\Delta'$. Subtracting the 
equality (\ref{DimEq}) from (\ref{DimEq2}) with $f=f'$ we get the reqiured equality of codimensions. ~~\bl

Now we apply this result to our computations of exceptional 2D fibers:

\subsubsection{Cells of the space ${\cal H}_2^2$.}\label{FibH22}    
Simplicial factor $\Delta[\Gamma]$ of the coordinate space
is a triangle when $[\Gamma]=[\pm\Gamma_2], [\Gamma_3]$. Periods map sends these triangles to the plane $(\Pi_0,\Pi_1)$ without degeneration so there are fibers of the maximal dimension in the boundary. However the boundary of the triangles is the projection of the exterior faces of the appropriate coordinate spaces.
The simplex $\Delta[\pm\Gamma_1]$ is three-dimensional and exactly one of its faces satisfies the condition of Lemma \ref{fDelta}.
Unfortunately, this face lifts to an outer face of ${\cal A}[\pm\Gamma_1]$. Summarizing, all fibers of the periods map transversally intersect the boundaries of  full dimensional cells of the labyrinth space ${\cal L}_2^2$.    

\subsubsection{Cells of the space ${\cal H}_2^1$.}\label{FibH21}
Simplicial factor $\Delta[\Gamma]$ of the coordinate space has no boundary simplexes that meet the condition of the Lemma \ref{fDelta} 
when $[\Gamma]=[\Gamma_1],~[\Gamma_2]$ (see left panel of Fig. \ref{PIH21}). In all the remaining cases such boundary simplexes exist but they correspond to the outer boundary of the coordinate space ${\cal A}[\Gamma]$, with two exceptions: 
$[\Gamma]=[\pm\Gamma_4]$. In each of the latter cases two codimension one boundary simplexes satisfy the condition of Lemma \ref{fDelta} and they lift to the inner faces of the coordinate space. In the case $[\Gamma]=[\Gamma_4]$ the exceptional values 
of the period map belong to interiors of two sides of the triangle $\Pi\Delta[\Gamma_4]$ incident to $\Pi^*=(0,0)$. For the case $[\Gamma]=[-\Gamma_4]$ the interior points in the sides of the triangle $\Pi\Delta[-\Gamma_4]$ incident to $\Pi^*=(-4\pi,4\pi)$
are also lifted to 2D fibers.

\section{Polyhedral model of fibers of Periods map}
We have seen that investigation of fibers of periods map requires different instruments at different zoom.
At the microscopic level we use differential geometry to find out that the fibers are smooth.
At the mesoscale the linear algebra and convex analysis say that the fibers within blocks commensurate with the size of the moduli space remain cells. At the cosmological scale -- in the labyrinth space --  we have to use combinatorics to analyse the   
global structure of the fibers composed of glued together polyhedra.

Decomposition of the labyrinth space $\cal L$ into polyhedra with the mapping $\Pi(\cdot)$ being a linear map 
in each of them gives the following natural receipt for the construction of the fiber of 
the global period mapping above its value $\Pi^*$.

\subsection{Algorithm}
STEP1: Given the value $\Pi^*$ and the topological invariants $g,k$ of real curve, list all braids $\beta\in Br_{g-k+1}$ satisfying the inclusions\footnote{For the $g=2$ case we consider in detail below the number of braids satisfying the inclusion is finite. However already 
for the $g=3$, $k=1$ case the number of suitable braids may be infinite, so one has to glue countably many polyhedra to 
obtain the model of the fiber}
$$
\beta\cdot\Pi^*\in\cup_\Gamma \Pi({\cal L}\hat{\cal A}[\Gamma]),  \qquad codim[\Gamma]=0. 
$$
Technically, it is convenient to consider the closure of the right hand side of the inclusion -- the images of the simplexes 
$\Delta[\Gamma]$ -- with subsequent seaving braids $\beta$ which send the value $\Pi^*$ to the image of the exterior faces of coordinate space.

STEP2:  For each $\beta$ and $\Gamma$ from the above step we find the section of the closed coordinate space 
$\hat{\cal A}[\Gamma]$ by the $g-$plane of equal values of the period map:
$$
\hat{\cal A}[\Gamma,\beta\cdot\Pi^*]=\{(H,W)\in\hat{\cal A}[\Gamma]:\qquad \Pi(H)=\beta\cdot\Pi^*\}. 
$$
For the values $\beta\cdot\Pi^*$ in the interior of the images of full dimensional dressed coordinate spaces we get the 
$g$- dimensional polyhedra; for the points $\beta\cdot\Pi^*$ in the boundary of the images we select local fibers  
of maximal dimension with the help of Lemma \ref{fDelta}. 

STEP3: Glue the arising $g$-dimensional polyhedra by the same equivalence relations we used in (\ref{LRep}) to construct the labyrinth space: 
$$
\Pi^{-1}(\Pi^*)=\mathop{\sqcup}\limits_{\beta,[\Gamma] }~\beta^{-1}\cdot{\cal L}\hat{\cal A}[\Gamma,\beta\cdot\Pi^*]/\sim, \qquad codim[\Gamma]=0.
$$

\subsection{First step for genus two curves}
For the number of ovals $k=2,3$, no braids appear and this step of the algorithm is trivial.
For the space ${\cal L}_2^1$ we suppose w.l.o.g. that the value $\Pi^*$ belongs to the image 
$\Delta:=a\cup b\cup c_\pm\cup d$ of the canonical lifts of all full dimensional cells.
Here $a,\dots,d$ are closed triangles shown on the left Fig. \ref{PIH21}.

\begin{lmm}
The orbit of a point $\Pi^*\in\Delta$ under the action of $Br_2$ has more than one point in the triangle $\Delta$  iff  
$\Pi^*\in c_\pm\cup d$. In the latter case the orbit has exactly one point in each of the triangles
$c_\pm\setminus d$.
\end{lmm}

{\bf Proof.} The action of the braid generator on the plane consists in translating the point $\Pi=(\Pi_1,\Pi_2)$
parallel to the line $\{\Pi_1+\Pi_2=0\}$ by the value of double euclidean distance to it:
$$
B(\Pi_1,\Pi_2)^t=(\Pi_1-(\Pi_1+\Pi_2),\Pi_2+(\Pi_1+\Pi_2))^t.
$$
Therefore the intersection of the triangle $\Delta$ with  a full orbit of Burau action of braids on a point $\Pi^*$ 
in our case has the appearance:
$$
\Pi, B\Pi, B^2\Pi, \dots, B^s\Pi, \qquad \Pi\in\Delta,\quad s\ge 0.
$$   
If $s>0$, then for the first point in the chain we have $B\Pi\in\Delta\not\ni B^{-1}\Pi$, which can be written as 
$$
\Pi\in B^{-1}\Delta\setminus B\Delta=c_+\setminus d.
$$
The last point in the chain satisfies the inclusion
$$
\beta^s\Pi\in B\Delta\setminus B^{-1}\Delta=c_-\setminus d,
$$
and all the rest points of the orbit lie in 
$$
 B^{-1}\Delta\cap B\Delta=d. ~~~~~\bl
$$

\subsection{Second step: carving the patches.}
\subsubsection{Space ${\cal L}_2^2$:} For  the value $\Pi^*$ in the interiors of each of three triangles $a_\pm ,b$ from the Fig. \ref{PIH22} we reconstruct the fiber of the periods map above it from Tab. \ref{Tab22}. On the interiors of intersections 
$a_\pm \cap b$ the only fibers are half-strips which belong to the spaces $\hat{\cal A}[\pm\Gamma_1]$ as it follows from Sect. \ref{FibH22}.

\subsubsection{Space ${\cal L}_2^1$:} For points $\beta\cdot\Pi^*$ in  the interiors of each of five triangles $a,b,c_\pm,d$ from the Fig. \ref{PIH21} the fibers of period map above them are listed in Tab. \ref{Tab21}. The analysis given in Sect. \ref{FibH21} shows that on the interfaces of the triangles the local fibers of periods map are listed in Tab \ref{FibInt21}.

\begin{table}
\begin{tabular}{c|c|c}
 Range of $\beta\cdot\Pi^*$& Polyhedron $\hat{\cal A}[\Gamma]$& Fiber inside the polyhedron\\
\hline
$Int(a\cap b)$ &$\hat{\cal A}[\Gamma_1]$ & rectangle\\
& $\hat{\cal A}[\pm\Gamma_5]$ & half-strips\\
$Int(b\cap c_\pm)$ &$\hat{\cal A}[\Gamma_1]$ & trapezoid\\
& $\hat{\cal A}[\Gamma_2]$, $\hat{\cal A}[\pm\Gamma_5]$ & half-strip\\
& $\hat{\cal A}[\pm\Gamma_4]$, $\hat{\cal A}[\pm\Gamma_6]$, $\partial_{int}\hat{\cal A}[\mp\Gamma_4]$ & sectors\\
$Int(c_\pm\cap d)$ &$\hat{\cal A}[\Gamma_1]$ & triangle\\
& $\hat{\cal A}[\Gamma_2]$ & half-strip\\
& $\partial_{int}\hat{\cal A}[\pm\Gamma_4]$ & sectors\\
$b\cap d$ &$\hat{\cal A}[\Gamma_1]$ & triangle\\
& $\hat{\cal A}[\Gamma_2]$, $\hat{\cal A}[-\Gamma_2]$ & half-strip\\
& $\partial_{int}\hat{\cal A}[\Gamma_4]$, $\partial_{int}\hat{\cal A}[-\Gamma_4]$ & sectors\\
\end{tabular}
\caption{\small  Local fibers of the periods map above interfaces of triangles $a,b,c_\pm,d$. }
\label{FibInt21}
\end{table}

\subsection{Third step: the patchwork} 
Having prepared the polygons $\hat{\cal A}[\Gamma, \beta\cdot\Pi^*]$ labeled by a graph and a braid, we assemble them together using 
the identification rule from (\ref{LRep}). For convenience, we first glue the polygons labeled by the same braid $\beta$. Typically, we identify the sides of different patches iff the sides  are labeled by the same graph $[\Gamma]$ of positive codimension. For the space ${\cal L}_2^3$ consisting of a unique coordinate space the job had been already done in Sect.  \ref{FiberH23}: the fiber is always an open quadrant.

\subsubsection{Assembly for fibers of the space ${\cal L}_2^2$}
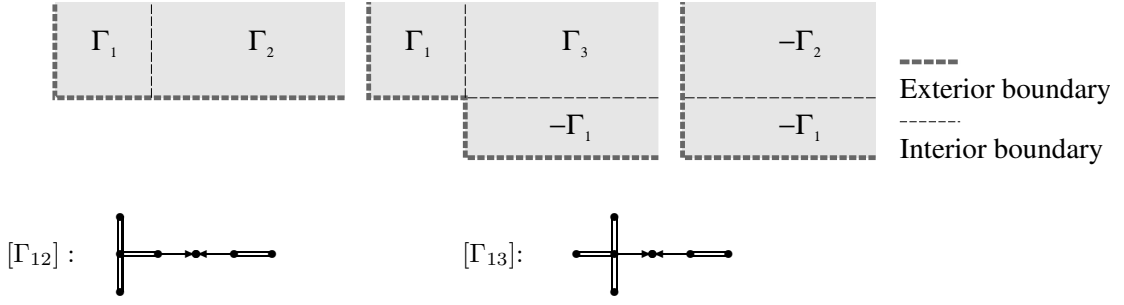
\begin{figure}[h!]
\centerline{\includegraphics[scale=.9, trim =-.4cm -.2cm 0cm 0cm, clip]{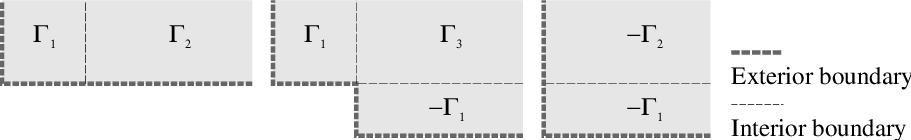}}
\begin{picture}(170,15)
\put(10,5){
\begin{picture}(45,20)
\multiput(5,0)(5,0){5}{\circle*{1}}
\multiput(5,0.3)(15,0){2}{\line(1,0){5}}
\multiput(5,-0.3)(15,0){2}{\line(1,0){5}}
\multiput(5,5)(0,-10){2}{\circle*{1}}

\put(4.7,5){\line(0,-1){10}}
\put(5.3,5){\line(0,-1){10}}

\put(15,0){\circle*{1}}
\put(10,0){\vector(1,0){5}}
\put(20,0){\vector(-1,0){5}}

\put(-10,-1){$[\Gamma_{12}]:$}
\end{picture}}

\put(70,5){
\begin{picture}(45, 15)
\multiput(5,0)(5,0){5}{\circle*{1}}
\multiput(5,0.3)(15,0){2}{\line(1,0){5}}
\multiput(5,-0.3)(15,0){2}{\line(1,0){5}}
\multiput(10,5)(0,-10){2}{\circle*{1}}
\put(15,0){\circle*{1}}

\put(9.7,-5){\line(0,1){10}}
\put(10.3,-5){\line(0,1){10}}

\put(10,0){\vector(1,0){5}}
\put(20,0){\vector(-1,0){5}}
\put(-10,-1){$[\Gamma_{13}]$:}
\end{picture}}
\end{picture}
\caption{\small  Intersection of the fiber $\Pi^{-1}(\Pi^*)$ with the full dimensional cells of the space ${\cal L}_2^2$.
Top row, left to right: $\Pi^*\in Int(a_+),~Int(b),~Int(a_-)$. The graphs on the interface of full dimensional cells are shown in the bottom row: $[\pm\Gamma_{12}]<[\pm\Gamma_1], [\Gamma_2]$;~~~ $[\pm\Gamma_{13}]<[\pm\Gamma_1], [\Gamma_3]$.} 
\label{FiberH22}
\end{figure}
Since the braids in case of genus two curves with two real ovals are trivial, the final answer is shown in Fig. \ref{FiberH22}
for the value  $\Pi^*$ in the interiors of the triangles $a_\pm,b$. On the interfaces of the triangles the global fiber of the periods map reduces to open half-strip.

\subsubsection{Assembly for fibers of the space ${\cal L}_2^1$}
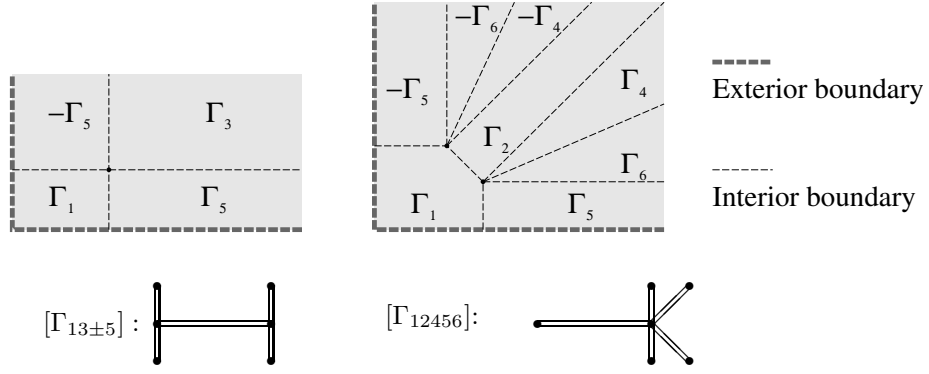
\begin{figure}[h!]
\centerline{\includegraphics[scale=.9, trim =-.4cm -.2cm 0cm 0cm, clip]{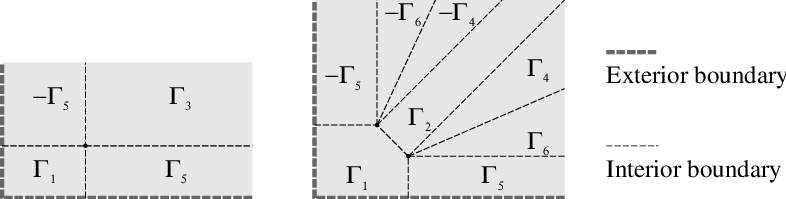}}
\begin{picture}(170,15)
\put(35,0){
\begin{picture}(45,20)
\multiput(0,0)(0,5){3}{\circle*{1}}
\multiput(15,0)(0,5){3}{\circle*{1}}
\multiput(-0.3,0)(0.6,0){2}{\line(0,1){10}}
\multiput(14.7,0)(0.6,0){2}{\line(0,1){10}}
\multiput(0,4.7)(0,0.6){2}{\line(1,0){15}}
\put(-15,4){$[\Gamma_{13\pm5}]:$}
\end{picture}}

\put(85,0){
\begin{picture}(45, 15)
\multiput(0,5)(15,0){2}{\circle*{1}}
\multiput(15,0)(5,0){2}{\circle*{1}}
\multiput(15,10)(5,0){2}{\circle*{1}}
\multiput(0,4.7)(0,.6){2}{\line(1,0){15}}
\multiput(14.7,0)(.6,0){2}{\line(0,1){10}}
\multiput(14.8,5.2)(.4,-.4){2}{\line(1,1){5}}
\multiput(14.8,4.8)(.4,.4){2}{\line(1,-1){5}}

\put(-20,5){$[\Gamma_{12456}]$:}
\end{picture}}
\end{picture}
\caption{\small  Top row: intersection of the fiber $\Pi^{-1}(\Pi^*)$ with the full dimensional cells of the space ${\cal L}_2^1$
when $\Pi^*\in Int~a$ (left) and $\Pi^*\in Int~b$ (right). ~~~~~~Bottom row:
Meeting point of several patches: $[\Gamma_{13\pm 5}]<[\Gamma_1], [\Gamma_3], [\Gamma_{\pm 5}]$ (left)
and $[\Gamma_{12456}]<[\Gamma_1], ~[\Gamma_2],\dots, ~[\Gamma_6]$ (right)} 
\label{FibH21ab}
\end{figure}

\begin{figure}[h!]
\centerline{\includegraphics[scale=.9, trim =-.4cm -.2cm 0cm 0cm, clip]{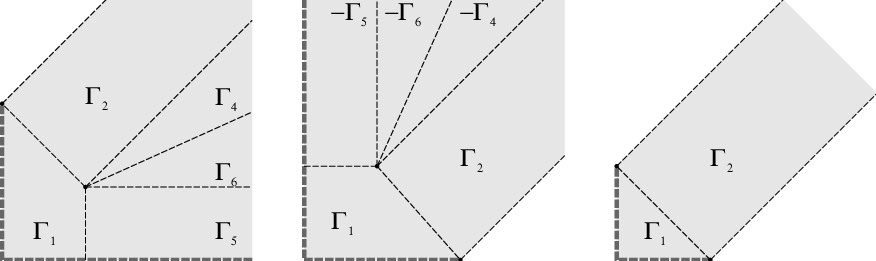}}
\caption{\small  Intersection of the fiber $\Pi^{-1}(\Pi^*)$ with the canonical lifts of full dimensional cells of the space ${\cal L}_2^1$
when $\Pi^*\in Int~ c_+$ (left), $Int~c_-$ (middle) or $Int~d$ (right).} 
\label{FibH21cde}
\end{figure}

First we assemble the patches  $\hat{\cal A}[\Gamma, \beta\cdot\Pi^*]$ with the same value of  $\beta\cdot\Pi^*$
in the interiors of the triangles $a,b,c_\pm,d$ shown in Fig. \ref{PIH21}: see Figs. \ref{FibH21ab}, \ref{FibH21cde}.

Fibers on the interfaces of the triangles are as follows:
\begin{itemize}
\item 
$Int(a\cap b)$, the same picture as on the left panel of Fig. \ref{FibH21ab} with removed patch corresponding 
to $\Gamma_3$ and totally exterior boundary. 
 \item 
$Int(b\cap c_\pm)$, same picture as for the fibers above triangle $Int(b)$
with two removed patches corresponding to $[\mp\Gamma_5]$, $[\mp\Gamma_6]$ and trapezoid block instead of the pentagonal
one corresponding to $[\Gamma_1]$. The patch corresponding to $[\pm\Gamma_4]$ lies on the boundary of the coordinate space and one side of the patch is the outer boundary.
 \item 
$Int(c_\pm\cap d)$  -- same pictures as for fibers above $Int(c_\pm)$
with two removed patches corresponding to $[\pm\Gamma_5]$, $[\pm\Gamma_6]$ and triangular block instead of the trapezoid
one corresponding to $[\Gamma_1]$. The patch corresponding to $[\pm\Gamma_4]$ lies on the boundary of the coordinate space and one side of this patch is the exterior boundary.
 \item 
$b\cap d$  -- the picture is the same as for fibers above $Int(b)$
with four removed patches corresponding to $[\pm\Gamma_5]$, $[\pm\Gamma_6]$ and triangular block instead of the pentagonal
one corresponding to $[\Gamma_1]$. The whole boundary of this patchwork is outer.
\end{itemize}

\begin{figure}[h!]
\centerline{\includegraphics[scale=.9, trim =0cm -.2cm 0cm 0cm, clip]{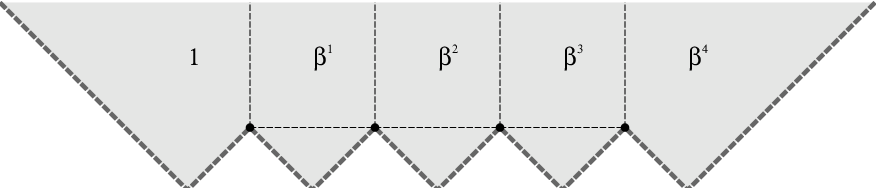}}
\caption{\small  Large patchwork, $\Pi^*\in c_+\cup c_-\cup d$.} 
\label{Glue}
\end{figure}

For the value $\Pi^*$ of the period map in $a\cup b$ the fiber is ready.
For the remaining cases $\Pi^*\in c_\pm\cup d$ one has to glue together pieces with different braid labels. 
The result of this operation is shown in Fig \ref{Glue}.

Summarizing, we see that any 2D fiber of the period map is a cell as it was stated in the main Theorem \ref{MainTh}.

\section{Conclusion}
We have developed an effective combinatorial approach to the study of the period mapping proposed in \cite{B2} and \cite{BBook}, Chaps. 4,5.
It had been proved that the global 2D fibers of the period mapping in the 4D labyrinth spaces have the simplest topology.
Basically, same technology may be applied for higher dimensional moduli spaces. The broad use of combinatorial techniques is 
vital for the investigation of already 3D fibers. Indeed, the number of full dimensional cells 
in the decomposition of 6D moduli space ${\cal H}_3^1$ is too large to consider fibers inside each cell with subsequent 
gluing them together. One has to use larger 'building blocks', say clusters of coordinate spaces with non-exceptional graphs $\Gamma$. 
As we mentioned already, a new effect appears in the 3D world: the fibers may become disconnected, however each component remains a cell.
We are going to describe the 3D calculations in a forthcoming publication.

\vspace{5mm}
\parbox{9cm}
{\it
119991 Russia, Moscow GSP-1, ul. Gubkina 8,\\
Institute for Numerical Mathematics,\\
Russian Academy of Sciences\\[2mm]
{\tt ab.bogatyrev@gmail.com}}

\end{document}